\documentclass[11pt] {article}    
\usepackage{amsmath}
\usepackage{amssymb}

\usepackage{theorem}
\numberwithin{equation}{section}

\title{Vanishing Viscosity Limits and Boundary Layers
for Circularly Symmetric 2D Flows}

\author{M.~C.~Lopes Filho\thanks{The first author was
supported in part by CNPq grant 302.102/2004-3. }\\
Departamento de Matematica, IMECC-UNICAMP\\
Caixa Postal 6065, Campinas, SP 13081-970, Brazil\\
mlopes@ime.unicamp.br 
\and A.~L.~Mazzucato\thanks{The second author was supported
in part by NSF grant DMS-0405803.}\\
Department of Mathematics,
Penn State University\\
McAllister Building\\
University Park, PA 16802, U.S.A.\\
alm24@psu.edu
\and H.~J.~Nussenzveig Lopes\thanks{The third author was
supported in part by CNPq grant 302.214/2004-6.}\\
Departamento de Matematica, IMECC-UNICAMP\\
Caixa Postal 6065, Campinas, SP 13081-970, Brazil \\
hlopes@ime.unicamp.br
\and Michael Taylor\thanks{The fourth author was supported
in part by NSF grant DMS-0456861.} \\
Department of Mathematics, University of North Carolina \\
CB \#3250, Phillips Hall\\
Chapel Hill, NC 27599, U.S.A.\\
met@math.unc.edu}

\date{}

\newcommand{\Dbar} {\overline{D}}
\newcommand{\Ombar} {\overline{\Omega}}
\newcommand{\ep} {\varepsilon}
\newcommand{\PI} {\operatorname{PI}}
\newcommand{\TV} {\operatorname{TV}}
\newcommand{\BV} {\operatorname{BV}}
\newcommand{\Tr} {\operatorname{Tr}}
\newcommand{\dv} {\operatorname{div}}
\newcommand{\rot} {\operatorname{rot}}

\newcommand{\loc} {\operatorname{loc}}
\newcommand{\pa} {\partial}
\newcommand{\RR} {\mathbb R}
\newcommand{\NN} {\mathbb N}

\newcommand{\Cal} {\mathcal}
\newcommand{\beq} {\begin{equation}}
\newcommand{\eeq} {\end{equation}}
\newcommand{\demo} {\noindent {\it Proof. }}
\newcommand{\qed} {\hfill$\Box$ {\newline $\text{}$}}

\theoremstyle{plain}
\newtheorem{theorem}{Theorem}[section]
\newtheorem{proposition}[theorem]{Proposition}

\usepackage{color}                                               %
\begin{document}

\maketitle

\begin{abstract}
We continue the work of Lopes Filho, Mazzucato and Nussenzveig Lopes
\cite{LMN} on the vanishing viscosity limit of circularly symmetric
viscous flow in a disk with rotating boundary, shown there to converge to
the inviscid limit in $L^2$-norm as long as the prescribed angular
velocity $\alpha(t)$ of the boundary has bounded total variation.
Here we establish convergence in stronger $L^2$ and $L^p$-Sobolev spaces,
allow for more singular angular velocities $\alpha$, and address the issue of
analyzing the behavior of the boundary layer.  This includes an analysis
of concentration of vorticity in the vanishing viscosity limit.
We also consider such flows on an annulus, whose two boundary components
rotate independently.  
\end{abstract}


\section{Introduction}\label{sec:1}

In this paper we study the 2D Navier-Stokes equation in the disk
$D=\{x\in\RR^2:|x|<1\}$:
\beq
\pa_t u^\nu+\nabla_{u^\nu}u^\nu+\nabla p^\nu=\nu\Delta u^\nu,\quad
\dv u^\nu=0,
\label{1.1}
\eeq
with no-slip boundary data on a rotating boundary:
\beq
u^\nu(t,x)=\frac{\alpha(t)}{2\pi} x^\perp,\quad |x|=1,\ t>0,
\label{1.2}
\eeq
and with circularly symmetric initial data:
\beq
u^\nu(0)=u_0(x),\quad \dv u_0=0,\quad u_0\, \| \, \pa D.
\label{1.3}
\eeq
In (\ref{1.2}), $x^\perp =Jx$, where $J$ is counterclockwise rotation by
$90^\circ$.  By definition, a vector field $u_0$ on $D$ is circularly
symmetric provided
\beq
u_0(R_\theta x)=R_\theta u_0(x),\quad \forall\, x\in D,
\label{1.4}
\eeq
for each $\theta\in [0,2\pi]$, where $R_\theta$ is counterclockwise rotation
by $\theta$.  The general vector field satisfying (1.4) has the form
\beq
s_0(|x|)x^\perp+s_1(|x|)x,
\label{1.5}
\eeq
with $s_j$ scalar,
but the condition $\dv u_0=0$, together with the condition $u_0\, \| \,\pa D$,
forces $s_1\equiv 0$, 
so the type of initial data we consider is characterized by
\beq
u_0(x)=s_0(|x|)x^\perp.
\label{1.6}
\eeq

Another characterization of vector fields of the form (\ref{1.6}) is the following.
For each unit vector $\omega\in S^1\subset\RR^2$, let
$\Phi_\omega:\RR^2\rightarrow \RR^2$
denote the reflection across the line generated by $\omega$, i.e.,
$\Phi_\omega(a\omega+bJ\omega)=a\omega-bJ\omega$.
Then a vector field $u_0$ on $D$ has the form (1.6)
if and only if
\beq
u_0(\Phi_\omega x)=-\Phi_\omega u_0(x),\quad \forall\, \omega\in S^1,\ x\in D.
\label{1.7}
\eeq

As is well known, such a vector field $u_0$ as given in (1.6) is a steady
solution to the 2D Euler equation.  In fact, a calculation gives
\beq
\nabla_{u_0}u_0=-s_0(|x|)^2x=-\nabla p_0(x),
\label{1.8}
\eeq
with
\beq
p_0(x)=\tilde{p}_0(|x|),\quad \tilde{p}_0(r)=-\int_r^1 \rho s_0(\rho)^2\,
d\rho,
\label{1.9}
\eeq
and the assertion follows.  We mention that $\|p_0\|_{L^1(D)}\le C\|u_0
\|^2_{L^2(D)}$.

The problem we address is the following {\it vanishing viscosity problem}:
to demonstrate that the solution $u^\nu$ to (\ref{1.1})--(\ref{1.3}) satisfies
\beq
\lim\limits_{\nu\searrow 0}\, u^\nu(t,\cdot)=u_0,
\label{1.10}
\eeq
and in particular to specify in what topologies such convergence holds.
As has been observed, what makes this problem tractable is the following
result.

\begin{proposition} \label{p1.1}
Given that $u_0$ has the form (\ref{1.6}), the solution
$u^\nu$ to (\ref{1.1})--(\ref{1.3}) is circularly symmetric for each $t>0$, of the
form
\beq
u^\nu(t,x)=s^\nu(t,|x|)x^\perp,
\label{1.11}
\eeq
and it coincides with the solution to the linear PDE
\beq
\pa_t u^\nu=\nu\Delta u^\nu,
\label{1.12}
\eeq
with boundary condition (\ref{1.2}) and initial condition (\ref{1.3}).
\end{proposition}

Here is a brief proof.  Let $u^\nu$ solve (\ref{1.12}), (\ref{1.2}), and (\ref{1.3}), with
$u_0$ as in (\ref{1.6}).  We claim (\ref{1.11}) holds.  In fact, for each unit vector
$\omega\in\RR^2$,
\beq
-\Phi_\omega u^\nu(t,\Phi_\omega x)
\label{1.13}
\eeq
also solves (\ref{1.12}) with the same initial data and boundary condition as
$u^\nu$, so these functions coincide, and (\ref{1.11}) follows.  Hence $\dv u^\nu
=0$ for each $t>0$.  Also we have an analogue of (\ref{1.8})--(\ref{1.9}):
\beq
\aligned
&\nabla_{u^\nu}u^\nu=-\nabla p^\nu,\quad p^\nu(t,x)=\tilde{p}^\nu(t,|x|), \\
&\tilde{p}^\nu(t,r)=-\int_r^1 \rho s^\nu(t,\rho)^2\, d\rho. 
\endaligned
\label{1.14}
\eeq
Hence this $u^\nu$ is the solution to (\ref{1.1})--(\ref{1.3}). 
For additional discussion of this issue, in particular in the context of weak
solutions, see \cite{LMN}, Proposition 5.1.


Previous work on the convergence problem (\ref{1.10}), in the circularly
symmetric context, was done by Matsui \cite{Mat}, who considered the case
$\alpha = 0$, without assuming compatibility of the initial velocity with
this boundary condition,
(see also \cite{Kell} for another treatment of the convergence problem in the 
circularly symmetric context), by Wang 
\cite{W}, whose general work on the convergence 
(building on results of Kato \cite{K}) is applicable to the
circularly symmetric case when $\alpha\in H^1_{\text{loc}}(\RR)$, by
three of us in \cite{LMN}, who treated
\beq
\alpha\in \BV(\RR)
\label{1.16}
\eeq
(supported in $\RR^+$),
and also by Bona and Wu \cite{BW}, who dealt with the special case
\beq
\alpha\equiv 0,\quad u_0\bigr|_{\pa D}=0.
\label{1.17}
\eeq
Results in these papers yield convergence in (\ref{1.10}), in $L^2(D)$-norm,
locally uniformly in $t$, when $u_0\in L^2(D)$ has the form (\ref{1.6}).
(In the special case (\ref{1.17}), convergence in stronger norms for more
regular $u_0$ was obtained in \cite{BW}; cf.~\S{10} of this paper for further
discussion of this case.)

In this paper we sharpen the treatment of the vanishing viscosity convergence
in several important respects.  For one, we go beyond $L^2$-norm convergence,
and establish norm convergence, under appropriate hypotheses, in $L^q$-Sobolev
spaces $H^{s,q}(D)$, when $sq<1$.  This is the maximal class of Sobolev
spaces for which such results could hold, since, without special compatibility
hypotheses such as (\ref{1.17}), the Sobolev space trace theorems forbid convergence
in higher norms.  The techniques we use to get such results also allow us to
treat driving motions $\alpha$ much more singular than in (\ref{1.16}); in fact,
we treat
\beq
\alpha\in L^{p'}(\RR),\quad p'\ge 1
\label{1.18}
\eeq
(supported in $\RR^+$).

In addition, we establish much stronger local convergence results, given
more regular data $u_0$.  On each compact subset $\Ombar$ of $D$, convergence
in (\ref{1.10}) holds in $H^k(\Omega)$ as long as $u_0$ is of class $H^k$ on a
neighborhood of $\Ombar$.  Furthermore, we give a precise analysis of the
boundary layer behavior of $u^\nu(t,x)$, as $\nu\searrow 0$, showing the
transitional behavior on a layer about $\pa D$ of thickness $\sim \nu^{1/2}$,
in case $u_0\in C^\infty(\Dbar)$, and more generally in case $u_0\in C(\Dbar)$.

It is a classical open problem whether solutions of the Navier-Stokes equations
in a bounded domain with no-slip boundary data converge to solutions of the Euler 
equations in the vanishing viscosity limit. The results obtained here may be regarded as 
an exploration of the difficulty involved in this problem by means of a nearly 
explicit example. In particular, we highlight two aspects of our results. 
First, we prove concentration of vorticity at the boundary. As is well-known, concentration  
of vorticity creates difficulties in treating the inviscid limit, see \cite{schochet}. 
Second, we obtain an expression of the total 
mass of vorticity present in the domain in terms of the angular acceleration of the boundary,
something which may be used as a sharp test for the accurate portrayal of the 
fluid-boundary interaction in high Reynold number numerical schemes.      

The structure of the rest of this paper is as follows.  In \S{\ref{sec:2}} we give a
general description of solutions to (\ref{1.12}), (\ref{1.2}), and (\ref{1.3}), 
with quite rough $\alpha$.  In \S\S{\ref{sec:3}--\ref{sec:5}} 
we establish convergence in (\ref{1.10}), first
in $L^2$-Sobolev spaces, then in $L^q$-Sobolev spaces for other values
of $q$, and then in certain Banach spaces of distributions (defined in \S{\ref{sec:5}}),
treating $\alpha$ of the form (\ref{1.18}), first for $p'>4$ (in \S{\ref{sec:3}}), then
for $p'>2$ (in \S{\ref{sec:4}}), and finally for all $p'\ge 1$ (in \S{\ref{sec:5}}).  
In \S{\ref{sec:6}} we digress to remark on the case when $\alpha$ is Brownian motion.

Section \ref{sec:7} treats strong convergence results away from $\pa D$.  In \S{\ref{sec:8}}
we produce estimates on the pressure $p^\nu$ appearing in (\ref{1.1}), making use
of the identities in (\ref{1.14}).  In \S{\ref{sec:9}} we examine the vorticity $\omega^\nu=
\rot u^\nu$, and contrast the local convergence to $\rot u_0$, on compact
subsets of $D$, with the global behavior.  In particular we analyze the
concentration of vorticity on $\pa D$ as $\nu\searrow 0$.  We devote \S{\ref{sec:10}}
to consideration of the special case (\ref{1.17}), and extend results of \cite{BW}.
In \S\S{\ref{sec:11}--\ref{sec:12}} we bring the theory of layer potentials to bear on 
the analysis of (\ref{1.12}), (\ref{1.2}), and (\ref{1.3}), and produce a sharp analysis of the
boundary layer behavior of $u^\nu(t,x)$, in case $u_0\in C^\infty(\Dbar)$, and more 
generally in case $u_0\in C(\Dbar)$.

In \S{\ref{sec:13}} we extend the scope of our investigations from the setting of the
disk $D$ to an annulus $\Cal{A}=\{x\in\RR^2:\rho<|x|<1\}$, for some
$\rho\in (0,1)$, allowing for independent rotations of the two components
of $\pa\Cal{A}$.  Thus the boundary condition (\ref{1.2}) is replaced by
\beq
\aligned
u^\nu(t,x)=\ &\frac{\alpha_1(t)}{2\pi} x^\perp,\quad |x|=1,\ t>0, \\
&\frac{\alpha_2(t)}{2\pi} x^\perp,\quad |x|=\rho,\ t>0.
\endaligned
\label{1.19}
\eeq
We establish an analogue of Proposition \ref{p1.1}, from which the extensions
of most of the results of \S\S{\ref{sec:2}--\ref{sec:12}} are straightforward,
though the extension of the material of \S{\ref{sec:9}} requires further work.

\section{Solutions with irregular driving motion $\alpha$}\label{sec:2}

As explained in the introduction, the analysis of the Navier-Stokes equation
in the circularly symmetric case is reduced to the analysis of the
initial-boundary problem
\beq
\pa_tu^\nu=\nu\Delta u^\nu,\quad u^\nu(0,x)=u_0(x),\quad u^\nu(t,x)=
\frac{\alpha(t)}{2\pi}x^\perp,\ \text{ for } x\in \pa D.
\label{2.1}
\eeq
In passing from (\ref{1.1})--(\ref{1.3}) to (\ref{2.1}), it was crucial to assume that $u_0$ 
had the form (\ref{1.6}), but such an hypothesis will not generally play an 
important role in our analysis of (\ref{2.1}), with some exceptions, such as in 
\S{\ref{sec:9}}.

We solve (\ref{2.1}) 
on $(t,x)\in \RR^+\times D$, but it is convenient to assume $\alpha$
is defined on $\RR$, with
\beq
\text{supp}\, \alpha\subset [0,\infty).
\label{2.2}
\eeq

As a preliminary to our main goal in this section of treating rough $\alpha$
in (\ref{2.1}), we first dispose of the case $\alpha\equiv 0$.  In this case,
if $u_0\in L^2(D)$, the solution to (\ref{2.1}) is given by
\beq
u^\nu(t)=e^{\nu tA}u_0,
\label{2.3}
\eeq
where $A$ is the self-adjoint operator on $L^2(D)$, with domain $\Cal{D}(A)$,
defined by
\beq
\Cal{D}(A)=H^2(D)\cap H^1_0(D),\quad Au=\Delta u\ \text{ for }\ u\in
\Cal{D}(A).
\label{2.4}
\eeq
Here we supress notation recording the fact that $u$ is vector-valued
(with values in $\RR^2$) rather than scalar valued.  

The family $\{e^{tA}:t\ge 0\}$ is a strongly continuous semigroup on $L^2(D)$.
As is well known, it also extends and/or restricts to a strongly continuous
semigroup on a large variety of other Banach spaces of functions on $D$,
such as $L^p(D)$ for $p\in [1,\infty)$ (but not for $p=\infty$).  The
maximum  principle holds; $\{e^{tA}:t\ge 0\}$ is a contraction semigroup
on $L^\infty(D)$, and also on $C(\Dbar)$, but these semigroups are not
strongly continuous at $t=0$.  We do get a strongly continuous semigroup
on $C_*(\Dbar)$, the space of functions in $C(\Dbar)$ that vanish on $\pa D$.
We also get strongly continuous semigroups on a variety of $L^q$-Sobolev
spaces, which we will discuss in more detail in \S\S{\ref{sec:3}--\ref{sec:4}}.  Whenever
$\{e^{tA}:t\ge 0\}$ acts as a strongly continuous semigroup on some Banach
space $X$ of functions on $D$, we get convergence in (\ref{2.3}) of $u^\nu(t)$
to $u_0$ in $X$-norm, for all $u_0\in X$.

In general, the solution to (\ref{2.1}) can be written as a sum of $e^{t\nu A}u_0$
and the restriction to $t\in [0,\infty)$ of a function that it is convenient
to define on $\RR\times D$ as the solution to
\beq
\aligned
&\pa_t v^\nu=\nu\Delta v^\nu,\quad v^\nu(t,x)=0\ \text{for}\ t<0, \\
&v^\nu(t,x)=\frac{\alpha(t)}{2\pi} x^\perp\ \text{for}\ x\in \pa D, \quad t\geq 0.
\endaligned
\label{2.5}
\eeq
Recall that we are assuming (\ref{2.2}).  We denote $v^\nu$ in (\ref{2.5}) as $\Cal{S}^\nu
\alpha$.  It is classical that
\beq
\Cal{S}^\nu:C^\infty_{\flat}(\RR)\longrightarrow C^\infty_{\flat} (\RR\times\Dbar),
\label{2.6}
\eeq
valid for each $\nu>0$.  Here and below, given a space $\mathfrak{X}$ of functions
or distributions on $\RR$ or $\RR\times D$, we denote by $\mathfrak{X}_{\flat}$ the
subspace consisting of elements of $\mathfrak{X}$ that vanish for $t<0$.  Thanks
to the maximum principle, $\Cal{S}^\nu$ in (\ref{2.6}) has a unique continuous
extension to
\beq
\Cal{S}^\nu:C_{\flat}(\RR)\longrightarrow C_{\flat}(\RR\times\Dbar).
\label{2.7}
\eeq
Our next goal is to show that $\Cal{S}^\nu$ also maps $L^{p'}_{\flat}(\RR)$
into other function spaces, on which the boundary trace $\Tr$ is defined,
and that
\beq
\Tr(\Cal{S}^\nu \alpha)=\frac{\alpha}{2\pi} x^\perp
\label{2.8}
\eeq
whenever $\alpha\in L^{p'}_{\flat}(\RR)$.

Note that in cases (\ref{2.6}) and (\ref{2.7}) $\Cal{S}^\nu\alpha$ clearly has a
boundary trace and (\ref{2.8}) holds.  Let us produce a variant of (\ref{2.7}) as
follows.  Using radial coordinates $(r,\theta)$ on $\Dbar$ (away from
the center), we have
\beq
\Cal{S}^\nu:C_{\flat}(\RR)\longrightarrow C([0,1],C_{\flat}(\RR\times\pa D)),
\label{2.9}
\eeq
i.e., $\Cal{S}^\nu\alpha(t,x)$, with $x=(r \cos \theta,r \sin \theta)$,
is a continuous function of $r\in [0,1]$ with values in the space
$C_{\flat}(\RR\times\pa D)$.  Then $\Tr\, \Cal{S}^\nu\alpha \in C_{\flat}(\RR\times
\pa D)$ is the value of this function at $r=1$.  Noting that $C_{\flat}(\RR
\times\pa D)\subset L^2_{\loc,\flat}(\RR\times\pa D)$, we have
\beq
\Cal{S}^\nu:C_{\flat}(\RR)\longrightarrow C([0,1],L^2_{\loc,\flat}(\RR\times\pa D)).
\label{2.10}
\eeq
Next note that if $\beta\in C_{\flat}^\infty(\RR)$ then
\beq
\alpha=\beta'\Longrightarrow \Cal{S}^\nu\alpha=\pa_t \Cal{S}^\nu\beta.
\label{2.11}
\eeq
From this it follows easily that
\beq
\Cal{S}^\nu\pa_t=\pa_t\Cal{S}^\nu:C^\infty_{\flat}(\RR)\longrightarrow
C([0,1],C^\infty_{\flat}(\RR\times\pa D))
\label{2.12}
\eeq
has a unique continuous extension to
\beq
\Cal{S}^\nu\pa_t=\pa_t\Cal{S}^\nu:C_{\flat}(\RR)\longrightarrow
C([0,1],H^{-1}_{\loc,\flat}(\RR\times\pa D)).
\label{2.13}
\eeq
Now, given $p'\ge 1$, each $\alpha\in L^{p'}_{\flat}(\RR)$ has the form
$\alpha=\beta'$ with $\beta\in C_{\flat}(\RR)$, namely $\beta(t)=\int_{-\infty}^t
\alpha(s)\, ds$.  It follows that
\beq
\Cal{S}^\nu:L^{p'}_{\flat}(\RR)\longrightarrow C([0,1],H^{-1}_{\loc,\flat}
(\RR\times\pa D)),
\label{2.14}
\eeq
for each $p'\ge 1$.  Consequently we have the continuous linear map
\beq
\Tr\circ \Cal{S}^\nu:L^{p'}_{\flat}(\RR)\longrightarrow
H^{-1}_{\loc,\flat}(\RR\times\pa D),
\label{2.15}
\eeq
and since (\ref{2.8}) holds on the dense linear subspace $C^\infty_0(\RR^+)$,
it holds for all $\alpha\in L^{p'}_{\flat}(\RR)$.

The target space in (\ref{2.14}) was chosen to have good trace properties, so
(2.8) could be verified, but such a choice precludes establishing the
convergence result
\beq
\lim\limits_{\nu\searrow 0}\, \Cal{S}^\nu \alpha=0
\label{2.16}
\eeq
in the strong topology of this target space.  Other spaces will arise
in \S\S{\ref{sec:3}--\ref{sec:5}}, for which (\ref{2.16}) holds in 
norm (see also (\ref{2.23})).  At this
point we will establish some useful identities for $\Cal{S}^\nu\alpha$.

To begin, we will assume $\alpha\in C^\infty_{\flat}(\RR)$; once we have the
identities we can extend the range of their validity by limiting arguments.
With $v^\nu$ defined by (\ref{2.5}), let us set
\beq
w^\nu(t,x)=v^\nu(t,x)-\frac{\alpha(t)}{2\pi} x^\perp
\label{2.17}
\eeq
on $[0,\infty)\times D$, so $w^\nu$ solves
\beq
\pa_t w^\nu=\nu\Delta w^\nu-\alpha'(t)f_1,\quad w^\nu(0,x)=0,\quad
w^\nu\bigr|_{\RR^+\times\pa D}=0,
\label{2.18}
\eeq
with
\beq
f_1(x)=\frac{1}{2\pi} x^\perp.
\label{2.19}
\eeq
We can then apply Duhamel's formula to write
\beq
w^\nu(t)=-\int_0^t e^{\nu(t-s)A}f_1\, \alpha'(s)\, ds.
\label{2.20}
\eeq
Hence
\beq
\aligned
\Cal{S}^\nu\alpha(t)=v^\nu(t)&=\alpha(t)f_1-\int_0^t e^{\nu(t-s)A}f_1\,
\alpha'(s)\, ds \\ &=\int_0^t \bigl(I-e^{\nu(t-s)A}\bigr)f_1\,
\alpha'(s)\, ds,
\endaligned
\label{2.21}
\eeq
and so the solution to (\ref{2.1}) can be written
\beq
u^\nu(t)=e^{\nu tA}u_0+\int_0^t \bigl(I-e^{\nu(t-s)A}\bigr)f_1\,
\alpha'(s)\, ds,
\label{2.22}
\eeq
with $f_1$ given by (\ref{2.19}).

Having (\ref{2.21}) and (\ref{2.22}) for $\alpha\in C^\infty_{\flat}(\RR)$, we can immediately
extend these formulas to $\alpha\in H^{1,1}_{\flat}(\RR)$, i.e., $\alpha$ supported
in $\RR^+$ and $\alpha,\ \alpha'\in L^1(\RR)$.  In fact, as in \cite{LMN}, we
can go further.  Let $X$ be a Banach space of functions on $D$ such that
$f_1\in X$ and $\{e^{tA}:t\ge 0\}$ is a strongly continuous semigroup on
$X$.  For example, we could take $X=L^2(D)$.  Here is the result.

\begin{proposition} \label{p2.1}
We have
\beq
\Cal{S}^\nu:\BV_{\flat}(\RR)\longrightarrow C_{\flat}(\RR,X),
\label{2.23}
\eeq
given by
\beq
\Cal{S}^\nu \alpha(t)=\int\limits_{I(t)} \bigl(I-e^{\nu(t-s)A}\bigr)f_1\,
d\alpha(s),\quad I(t)=[0,t],
\label{2.24}
\eeq
the integral being the Bochner integral.  
\end{proposition}
\demo
Using mollifiers in $C_0^\infty(0,1/k)$, we can approximate $\alpha$ by
$\alpha_k\in C^\infty_{\flat}(\RR)$ in a fashion so that $\alpha'\rightarrow d\alpha$
$\text{weak}^*$ as Radon measures.  That is to say,
\beq
\lim\limits_{k\rightarrow \infty}\, \int\limits_{\RR^+} g(s)\alpha_k'(s)\, ds
=\int\limits_{\RR^+} g(s)\, d\alpha(s),\quad \RR^+=[0,\infty),
\label{2.24A}
\eeq
for each compactly supported continuous function $g$ on $[0,\infty)$.
We have seen that $\Cal{S}^\nu\alpha_k\rightarrow \Cal{S}^\nu\alpha$, and that
$\Cal{S}^\nu\alpha_k$ is given by (\ref{2.22}) with $\alpha$ replaced by
$\alpha_k$.  To prove (\ref{2.24}), it suffices to show that
\beq
\aligned
\lim\limits_{k\rightarrow \infty}\, &\int_0^t \bigl\langle 
\bigl(I-e^{\nu(t-s)A}\bigr)f_1,\xi\bigr\rangle \alpha'_k(s)\, ds \\
&=\int\limits_{I(t)} \bigl\langle\bigl(I-e^{\nu(t-s)A}\bigr)f_1,\xi
\bigr\rangle\, d\alpha(s),
\endaligned
\label{2.24B}
\eeq
for arbitrary $\xi\in X'$.  However, (\ref{2.24B}) follows from (\ref{2.24A})
upon taking
\beq
\aligned
g(s)=\bigl\langle\bigl(I-e^{\nu(t-s)A}\bigr)f_1,\xi\bigr\rangle\quad
&\text{for }\ s\in [0,t], \\
0\qquad {}\qquad {}\quad {} &\text{for }\ s>t,
\endaligned
\label{2.24C}
\eeq
which is continuous and compactly supported on $[0,\infty)$.
Finally, the fact that the range in (\ref{2.23}) is contained in $C_{\flat}(\RR,X)$
is a consequence of the formula (\ref{2.24}).
\qed

\noindent{\it Remark.}
Note that the continuous integrand in (\ref{2.24})
vanishes at $s=t$, so one gets the same result with $I(t)=[0,t)$.
Note also that
\beq
\|\Cal{S}^\nu\alpha(t)\|_X\le \|\alpha\|_{\BV([0,t])}\,
\sup\limits_{s\in [0,t]}\, \|e^{\nu sA}f_1-f_1\|_X,
\label{2.25}
\eeq
and, if $u_0\in X$,
\beq
\|u^\nu(t)-u_0\|_X\le \|e^{\nu tA}u_0-u_0\|_X+\|\Cal{S}^\nu\alpha(t)\|_X.
\label{2.26}
\eeq
With $X=L^2(D)$, this gives the convergence result established in \cite{LMN}.

$\text{}$

Another useful identity for $\Cal{S}^\nu\alpha$ arises via integration by
parts.  In fact, for $\alpha\in C^\infty_{\flat}(\RR)$ and $\ep>0$, we have
\beq
\aligned
&\int_0^{t-\ep}e^{\nu(t-s)A}f_1\, \alpha'(s)\, ds \\
&=\alpha(t-\ep)e^{\nu\ep A}f_1+\nu\int_0^{t-\ep} Ae^{\nu(t-s)A}f_1\,
\alpha(s)\, ds,
\endaligned
\label{2.27}
\eeq
justification following because $e^{\nu\sigma A}f_1\in \Cal{D}(A)$, given
by (\ref{2.4}), whenever $\sigma>0$.  Together with (\ref{2.21}), this yields
\beq
\Cal{S}^\nu\alpha(t)=-\lim\limits_{\ep\searrow 0}\, \nu \int_0^{t-\ep}
Ae^{\nu(t-s)A}f_1\, \alpha(s)\, ds,
\label{2.28}
\eeq
the limit existing certainly in $L^2$-norm, locally uniformly in $t$,
and as we will see in subsequent sections, also in other norms, and for
more singular $\alpha$.

\section{$L^2$-Sobolev vanishing viscosity limits}\label{sec:3}

The family $\{e^{tA}:t\ge 0\}$ is a strongly continuous semigroup of
operators on $L^2(D)$, and on $\Cal{D}(A)$, given by (\ref{2.4}), and more generally
on $\Cal{D}((-A)^{\sigma/2})$, for each $\sigma\in\RR^+$.  As is well known,
\beq
\Cal{D}((-A)^{1/2})=H^1_0(D),
\label{3.1}
\eeq
and, for $\sigma\in [0,1]$,
\beq
\Cal{D}((-A)^{\sigma/2})=[L^2(D),H^1_0(D)]_\sigma,
\label{3.2}
\eeq
the complex interpolation space.  Furthermore,
\beq
\aligned
&[L^2(D),H^1_0(D)]_\sigma &=\ H^\sigma_0(D),\quad \frac{1}{2}<\sigma\le 1, \\
&\, &=\  H^\sigma(D),\quad 0\le\sigma<\frac{1}{2}.
\endaligned
\label{3.3}
\eeq
Cf.~\cite{LM}.  Consequently,
\beq
\Cal{D}((-A)^{\sigma/2})=H^\sigma(D),\quad \text{for }\ \sigma\in\Bigl[
0,\frac{1}{2}\Bigr).
\label{3.4}
\eeq
Hence
\beq
\forall\, \sigma\in\Bigl[0,\frac{1}{2}\Bigr),\quad u_0\in H^\sigma(D)
\Longrightarrow e^{\nu tA}u_0\rightarrow u_0\ \text{in}\
H^\sigma\text{-norm, as}\ \nu\rightarrow 0,
\label{3.5}
\eeq
convergence holding uniformly in $t\in [0,T]$ for each $T<\infty$.

Recall the formula (\ref{2.28}), i.e.,
\beq
\Cal{S}^\nu\alpha(t)=-\lim\limits_{\ep\searrow 0}\, \nu \int_0^{t-\ep}
Ae^{\nu(t-s)A}f_1\, \alpha(s)\, ds,
\label{3.6}
\eeq
for $\alpha\in C^\infty_{\flat}(\RR)$.  We now seek conditions that imply
that $\|A e^{\nu(t-s)A}f_1 \alpha(s)\|_{H^\sigma}$ is integrable on
$[0,t]$ and that
\beq
\Cal{S}^\nu\alpha(t)=-\nu \int_0^t Ae^{\nu(t-s)A}f_1\, \alpha(s)\, ds.
\label{3.7}
\eeq
Note that
\beq
\aligned
&\int_0^t \|\nu Ae^{\nu(t-s)A}f_1\,\alpha(s)\|_{H^\sigma(D)}\, ds \\
&\le \|\alpha\|_{L^{p'}([0,t])} \Bigl(\int_0^t
\|\nu A e^{\nu sA}f_1\|^p_{H^\sigma(D)}\, ds\Bigr)^{1/p}.
\endaligned
\label{3.8}
\eeq
Now
\beq
\aligned
&\|\nu Ae^{\nu sA}f_1\|_{H^\sigma(D)} \\
&\le C\|\nu(-A)^{1+\sigma/2}e^{\nu sA}f_1\|_{L^2(D)} \\ 
&=C\|\nu(-A)^{1-(\tau-\sigma)/2}e^{\nu sA}(-A)^{\tau/2}f_1\|_{L^2(D)} \\ 
&=C\nu^{(\tau-\sigma)/2} s^{(\tau-\sigma)/2-1}
\|(-\nu sA)^{1-(\tau-\sigma)/2} e^{\nu sA} (-A)^{\tau/2}f_1\|_{L^2(D)} \\
&\le C\nu^{(\tau-\sigma)/2}s^{(\tau-\sigma)/2-1} \|f_1\|_{H^\tau(D)}.
\endaligned
\label{3.9}
\eeq
Furthermore,
\beq
\aligned
&p\Bigl(1-\frac{\tau-\sigma}{2}\Bigr)<1 \Longrightarrow \\
&\int_0^t (s^{(\tau-\sigma)/2-1})^p\, ds=C_{p\sigma\tau}
t^{p(\tau-\sigma)/2-p+1},\quad C_{p\sigma\tau}<\infty.
\endaligned
\label{3.10}
\eeq
Hence
\beq
\aligned
&1\le p<\frac{2}{2-(\tau-\sigma)} \\
&\Longrightarrow
\|\Cal{S}^\nu\alpha(t)\|_{H^\sigma(D)}\le C(t) \nu^{(\tau-\sigma)/2}
\|\alpha\|_{L^{p'}([0,t])} \|f_1\|_{H^\tau(D)}.
\endaligned
\label{3.11}
\eeq
Approximating a rough $\alpha$ by smoothing convolutions and passing to
the limit, we obtain:

\begin{proposition} \label{p3.1} 
Assume $\alpha\in L^{p'}_{\flat}(\RR)$
where $p$ satisfies the hypothesis in (\ref{3.11}). 
Then the formula (\ref{3.7}) holds, we have
\beq
\aligned
&\Cal{S}^\nu:L^{p'}_{\flat}(\RR)\longrightarrow C_{\flat}(\RR,H^\sigma(D)), \\
&\text{for }\ \sigma\in \Bigl[0,\frac{1}{2}\Bigr),\ p\in \Bigl[1,
\frac{2}{3/2+\sigma}\Bigr),
\endaligned
\label{3.12}
\eeq
and the estimate (\ref{3.11}) holds.
\end{proposition}

$\text{}$ \newline
{\it Remark.}  There exist $\sigma,\tau$ such that $0\le\sigma<\tau<1/2$ and
the hypothesis above on $p$ holds provided $1\le p<4/3$, i.e., provided
$p'>4$.

\section{$L^q$-Sobolev vanishing viscosity limits}\label{sec:4}

Let $q\in (1,\infty)$.  Then $e^{tA}$ provides a strongly continuous
semigroup on $L^q(D)$, indeed a holomorphic semigroup.  We sometimes
denote the infinitesimal generator by $A_q$, to emphasize the $q$-dependence.
Now, for $\lambda>0$,
\beq
R_\lambda=(\lambda-A)^{-1}=\int_0^\infty e^{tA}e^{-\lambda t}\, dt
\label{4.1}
\eeq
has the mapping property
\beq
R_\lambda:L^q(D)\overset{\approx}{\longrightarrow}
\Cal{D}(A_q),
\label{4.2}
\eeq
and standard elliptic theory gives
\beq
\Cal{D}(A_q)=H^{2,q}(D)\cap H^{1,q}_0(D).
\label{4.3}
\eeq
We record the following useful known results.

\begin{proposition} \label{p4.1}
Given $\sigma\in (0,2)$, the operator
$-(-A_q)^{\sigma/2}$ is well defined and is the generator of a holomorphic
semigroup on $L^q(D)$.  Furthermore,
\beq
\Cal{D}((-A_q)^{\sigma/2})=[L^q(D),\Cal{D}(A_q)]_{\sigma/2},
\label{4.4}
\eeq
where the right side is a complex interpolation space.  In addition,
\beq
0\le\sigma<\frac{1}{q}\Longrightarrow \Cal{D}((-A_q)^{\sigma/2})=
H^{\sigma,q}(D).
\label{4.5}
\eeq
Also, if $\gamma\in [0,1]$ and $T\in (0,\infty)$,
\beq
\|(-tA_q)^{\gamma} e^{tA_q}f\|_{L^q(D)}\le
C_{q\gamma T}\, \|f\|_{L^q(D)},\quad \text{for }\ t\in [0,T].
\label{4.6}
\eeq
\end{proposition}
\demo
The results (\ref{4.4})--(\ref{4.5}) are proven in \cite{Se1}--\cite{Se2}.  
The result (\ref{4.6}) is equivalent to
\beq
\|e^{tA_q}f\|_{\Cal{D}((-A_q)^\gamma)}\le C t^{-\gamma}\|f\|_{L^q(D)}.
\label{4.7}
\eeq
For $\gamma=1$ this follows from the fact that $e^{tA_q}$ is a holomorphic
semigroup.  For $\gamma=0$ it is clear.  Then for $0<\gamma<1$ it follows 
from these endpoint cases, via (\ref{4.4}) and the general interpolation estimate
\beq
\|g\|_{[L^q,\Cal{D}(A_q)]_\gamma}\le C\|g\|^\gamma_{\Cal{D}(A_q)}
\|g\|_{L^q(D)}^{1-\gamma}.
\label{4.8}
\eeq
\qed

Given this proposition, we proceed as follows on the estimation of
\beq
\Cal{S}^\nu\alpha(t)=-\int_0^t \nu Ae^{\nu(t-s)A}f_1\, \alpha(s)\, ds.
\label{4.9}
\eeq
Pick $q\in (1,\infty)$, and pick $\sigma,\tau$ satisfying
\beq
0\le\sigma<\tau<\frac{1}{q}.
\label{4.10}
\eeq
Then, as in (\ref{3.8}), we have
\beq
\|\Cal{S}^\nu\alpha(t)\|_{H^{\sigma,q}(D)}\le \|\alpha\|_{L^{p'}([0,t])}
\Bigl(\int_0^t \|\nu Ae^{\nu sA}f_1\|^p_{H^{\sigma,q}(D)}\, ds\Bigr)^{1/p}.
\label{4.11}
\eeq
Now Proposition \ref{p4.1} yields the following analogue of (\ref{3.9}):
\beq
\aligned
&\|\nu Ae^{\nu sA}f_1\|_{H^{\sigma,q}(D)} \\
&=C\|\nu(-A)^{1+\sigma/2}e^{\nu sA}f_1\|_{L^q(D)} \\
&=\|\nu(-A)^{1-(\tau-\sigma)/2}e^{\nu sA}(-A)^{\tau/2}f_1\|_{L^q(D)} \\
&=C\nu^{(\tau-\sigma)/2} s^{(\tau-\sigma)/2-1}
\|(-\nu sA)^{1-(\tau-\sigma)/2} e^{\nu sA}(-A)^{\tau/2} f_1\|_{L^q(D)} \\&
\le C\nu^{(\tau-\sigma)/2} s^{(\tau-\sigma)/2-1} \|f_1\|_{H^{\tau,q}(D)}.
\endaligned
\label{4.12}
\eeq

We can then use (\ref{3.10}) to conclude:

\begin{proposition} \label{p4.2}
We have
\beq
\aligned
&\Cal{S}^\nu:L^{p'}_{\flat}(\RR)\longrightarrow C_{\flat}(\RR,H^{\sigma,q}(D)), \\
&\text{for }\ q>1,\ \sigma\in\Bigl[0,\frac{1}{q}\Bigr),\
p\in\Bigl[1,\frac{2}{2-1/q+\sigma}\Bigr),
\endaligned
\label{4.13}
\eeq
and as long as (\ref{4.10}) holds,
\beq
\aligned
&1\le p<\frac{2}{2-(\tau-\sigma)} \\
&\Longrightarrow
\|\Cal{S}^\nu\alpha(t)\|_{H^{\sigma,q}(D)}
\le C(t)\nu^{(\tau-\sigma)/2} \|\alpha\|_{L^{p'}([0,t])}
\|f_1\|_{H^{\tau,q}(D)}.
\endaligned
\label{4.14}
\eeq
\end{proposition}

$\text{}$ \newline {\it Remark.}
Note that, for a given $p$, there exist $q,\ \tau$, and $\sigma$
satisfying (\ref{4.10}), for which the hypothesis in (\ref{4.14}) holds, 
provided $1\le p<2$, i.e., provided $p'>2$.

$\text{}$

In the setting of Proposition \ref{p4.1}, $e^{tA_q}$ is also a strongly continuous 
semigroup on $\Cal{D}((-A_q)^{\sigma/2})$, hence in the setting of (\ref{4.5}),
on $H^{\sigma,q}(D)$, so we also have:

\begin{proposition} \label{p4.3}
If $q\in (1,\infty)$ and $\sigma\in [0,1/q)$,
then
\beq
u_0\in H^{\sigma,q}(D)\Longrightarrow \lim\limits_{\nu\searrow 0}\
\|e^{\nu tA}u_0-u_0\|_{H^{\sigma,q}(D)}=0.
\label{4.15}
\eeq
\end{proposition}

\section{Generalized function space vanishing viscosity limits}\label{sec:5}

Here we show that for $\alpha\in L^{p'}_{\flat}(\RR)$, we have
$\Cal{S}^\nu\alpha(t)\rightarrow 0$ in various topologies weaker than the
$L^2$-norm, even when $p'\in [1,2]$.  We will use a continuation of the
scale $\Cal{D}((-A_2)^{s/2})=\Cal{D}_s$.  There are analogous results
involving $A_q$, which we will not discuss here.  As stated before, we have
\beq
\Cal{D}_2=H^2(D)\cap H^1_0(D).
\label{5.1}
\eeq
Also
\beq
0\le s\le 2\Longrightarrow \Cal{D}_s=[L^2(D),\Cal{D}_2]_{s/2},
\label{5.2}
\eeq
and in particular
\beq
0\le s<\frac{1}{2}\Longrightarrow \Cal{D}_s=H^s(D).
\label{5.3}
\eeq
For $s<0$, we set
\beq
\Cal{D}_s=\Cal{D}^*_{-s}.
\label{5.4}
\eeq
Details on this are given in Chapter 5, Appendix A of \cite{T}.  We mention that
\beq
\Cal{D}_{-1}=H^{-1}(D).
\label{5.5}
\eeq
Also
\beq
s=-\sigma<0\Longrightarrow \Cal{D}_s=(-A_2)^{\sigma/2} L^2(D).
\label{5.6}
\eeq

Given this, we have in parallel with (\ref{3.8})--(\ref{3.9}) that for $\sigma\in\RR$,
\beq
\|\Cal{S}^\nu\alpha(t)\|_{\Cal{D}_\sigma}\le \|\alpha\|_{L^{p'}([0,t])}
\Bigl(\int_0^t \|\nu A e^{\nu sA}f_1\|^p_{\Cal{D}_\sigma}\, ds\Bigr)^{1/p},
\label{5.7}
\eeq
for $p\in (1,\infty)$, and
\beq
\|\Cal{S}^\nu\alpha(t)\|_{\Cal{D}_\sigma}\le \|\alpha\|_{L^1([0,t])}\,
\sup\limits_{0\le s\le t}\, \|\nu Ae^{\nu sA}f_1\|_{\Cal{D}_\sigma},
\label{5.8}
\eeq
and for $-\infty<\sigma<\tau<1/2$,
\beq
\|\nu A e^{\nu sA}f_1\|_{\Cal{D}_\sigma}\le C\nu^{(\tau-\sigma)/2}
s^{(\tau-\sigma)/2-1} \|f_1\|_{\Cal{D}_\tau}.
\label{5.9}
\eeq
We still have (\ref{3.10}), and now we can take $\sigma<0$, as well as $\tau$ close
to $1/2$.  In particular, we have
\beq
-2<\sigma<-\frac{3}{2},\ \tau=\sigma+2
\Longrightarrow \|\Cal{S}^\nu\alpha(t)\|_{\Cal{D}_\sigma}
\le C\nu \|\alpha\|_{L^1([0,t])}\, \|f_1\|_{H^\tau(D)}.
\label{5.10}
\eeq

\section{A stochastic interlude}\label{sec:6}

Here, instead of having $\alpha$ be deterministic, we consider
$$
\alpha(s)=\omega(s),
$$
where $\omega\in\Cal{P}_0$, the space of continuous paths from $[0,\infty)$
to $\RR$ (such that $\omega(0)=0$) endowed with Wiener measure $W_0$, and
expectation $E_0$.

The estimates of \S{\ref{sec:3}} apply to $\Cal{S}^\nu\omega(t)$,
but we record special results for this stochastic situation.

We are dealing with
\beq
\Cal{S}^\nu(t,\omega)=\int_0^t (I-e^{\nu(t-s)A})f_1\, d\omega(s),
\label{6.1}
\eeq
which is a Wiener-Ito integral.  The integrand is independent of $\omega$,
so the analysis of such an integral is relatively elementary.  We make the
following:

$\text{}$ \newline
{\it Hypothesis 6.1} $H$ is a Hilbert space of functions on $D$, with
values in $\RR^2$, such that $f_1\in H$ and $\{e^{sA}:s\ge 0\}$ is a strongly
continuous semigroup on $H$.
\newline $\text{}$

In such a case, we have
\beq
E_0\Bigl(\|\Cal{S}^\nu(t,\cdot)\|^2_H\Bigr)=\int_0^t \|(e^{\nu sA}-I)f_1\|_H^2
\, ds.
\label{6.2}
\eeq
This is a standard identity in the scalar case (cf.~\cite{T}, Chapter 11,
Proposition 7.1, where an extra factor arises due to an idiosyncratic
normalization of Wiener measure made there), and extends readily to
Hilbert space-valued integrands, by taking an orthonormal basis of $H$
and examining each component.  As seen in \S{\ref{sec:3}}, this is applicable for
\beq
H=H^\tau(D,\RR^2),\quad 0\le\tau<\frac{1}{2}.
\label{6.3}
\eeq

It might be interesting to obtain statistical information on the boundary
layers that arise for $\Cal{S}^\nu(t,\omega)$.

\section{Local convergence results}\label{sec:7}

Here we examine convergence of
\beq
u^\nu(t)=e^{\nu tA}u_0+\Cal{S}^\nu\alpha(t)
\label{7.1}
\eeq
as $\nu\rightarrow 0$ to $u_0$, on compact subsets of $D$.  Recall that
\beq
\Cal{S}^\nu\alpha(t)= - \nu \int_0^t Ae^{\nu(t-s)A}f_1\, \alpha(s)\, ds,
\label{7.2}
\eeq
where $f_1$ is given by (\ref{2.19}), so
\beq
f_1\in C^\infty(\Dbar).
\label{7.3}
\eeq

We will prove the following.

\begin{proposition} \label{p7.1}
Assume that
\beq
u_0\in L^2(D),\quad u_0\bigr|_{\Cal{O}}\in H^k(\Cal{O}),\quad
\alpha\in L^1_{\flat}(\RR),
\label{7.4}
\eeq
where $\Cal{O}$ is an open subset of $D$, and take $\Ombar\subset\subset
\Cal{O}$.  Then
\beq
\lim\limits_{\nu\searrow 0}\, u^\nu(t)\bigr|_\Omega=u_0\bigr|_\Omega
\quad \text{in }\ H^k(\Omega),
\label{7.5}
\eeq
uniformly for $t\in [0,T]$, given $T<\infty$.
\end{proposition}
\demo
First we show
\beq
\lim\limits_{\nu\searrow 0}\, e^{\nu tA}u_0\bigr|_\Omega=u_0\bigr|_\Omega
\quad \text{in }\ H^k(\Omega).
\label{7.6}
\eeq
To see this, take $\varphi\in C_0^\infty(\Cal{O})$ such that $\varphi=1$ on a
neighborhood $\Omega_1$ of $\Ombar$, and write
\beq
u_0=\varphi u_0+(1-\varphi)u_0=u_1+u_2,
\label{7.7}
\eeq
so
\beq
u_1\in H^k_0(D)\subset \Cal{D}((-A)^{k/2}),\quad
u_2\in L^2(D),\quad u_2=0\ \text{ on }\ \Omega_1.
\label{7.8}
\eeq
It follows that
\beq
\lim\limits_{\nu\searrow 0}\, e^{\nu tA}u_1=u_1\ \text{ in }\
\Cal{D}((-A)^{k/2})\subset H^k(D),
\label{7.9}
\eeq
so we will have (\ref{7.6}) if we show that
\beq
\lim\limits_{\nu\searrow 0}\, e^{\nu tA}u_2\bigr|_{\Ombar}=0
\ \text{ in }\ C^\infty(\Ombar).
\label{7.10}
\eeq
To do this, we define $w(s,x)$ on $\RR\times\Omega_1$ by
\beq
\aligned
w(s,x)=e^{sA}u_2(x),\quad &s\ge 0, \\
0,\quad {}\quad &s<0.
\endaligned
\label{7.11}
\eeq
Then $w$ is a weak solution of $(\pa_s-\Delta)w=0$ on $\RR\times\Omega_1$,
and the well known hypoellipticity of $\pa_s-\Delta$ implies
\beq
w\in C^\infty(\RR\times\Omega_1).
\label{7.12}
\eeq
This implies (\ref{7.10}) and hence we have (\ref{7.6}).

Now that, once we have (\ref{7.6}), we can apply this to $f_1$ in place of $u_0$
and deduce that
\beq
\lim\limits_{\nu\searrow 0}\, e^{\nu(t-s)A}f_1\bigr|_{\Omega}=
f_1\bigr|_{\Omega}\ \text{ in }\ H^{k+2}(\Omega),
\label{7.13}
\eeq
uniformly on $0\le s\le t\le T$, which by (\ref{7.2}) then gives
\beq
\lim\limits_{\nu\searrow 0}\, \Cal{S}^\nu\alpha(t)\bigr|_\Omega=0\
\text{ in }\ H^k(\Omega),
\label{7.14}
\eeq
and hence proves (\ref{7.5}).
\qed

\section{Pressure estimates}\label{sec:8}

For $u^\nu$ given by (\ref{1.1})--(\ref{1.4}), the pressure gradient $\nabla p^\nu$
is given by
\beq
\nabla p^\nu=-\nabla_{u^\nu}u^\nu.
\label{8.1}
\eeq
It is convenient to rewrite the right side of (\ref{8.1}), using the general
identity
\beq
\nabla_v u=\text{div}(u\otimes v)-(\text{div}\, v)u
\label{8.2}
\eeq
(cf.~\cite{T}, Chapter 17, (2.43)), which in the current context yields
\beq
\nabla p^\nu=-\text{div}\, (u^\nu\otimes u^\nu).
\label{8.3}
\eeq

Recall that
\beq
u^\nu=e^{\nu tA}u_0+\int_0^t \bigl(I-e^{\nu(t-s)A}\bigr)f_1\, d\alpha(s),
\label{8.4}
\eeq
where
\beq
f_1(x)=\frac{1}{2\pi}x^\perp.
\label{8.5}
\eeq
For simplicity we will work under the assumption that $\alpha$ has bounded
variation on each interval $[0,t]$.  We will assume
\beq
u_0\in L^\infty(D)\cap H^{\tau,q}(D),
\label{8.6}
\eeq
with
\beq
q\in (1,\infty),\quad 0<\tau<\frac{1}{q}.
\label{8.7}
\eeq
We aim to prove the following.

\begin{proposition} \label{p8.1}
Let $u^\nu$ be given by (\ref{1.1})--(\ref{1.4}), and assume
$u_0$ satisfies (\ref{8.6})--(\ref{8.7}).  Assume $\alpha$ has locally bounded
variation on $[0,\infty)$.  Take $T\in (0,\infty)$ and $\nu_0>0$.  Then,
uniformly for $t\in [0,T]$, we have
\beq
u^\nu(t)\otimes u^\nu(t)\ \text{ bounded in }\ L^\infty(D)\cap
H^{\tau,q}(D),
\label{8.8}
\eeq
for $\nu\in (0,\nu_0]$, and, as $\nu\rightarrow 0$,
\beq
u^\nu(t)\otimes u^\nu(t)\longrightarrow u_0\otimes u_0\ 
\text{ weak}^* \text{ in }\
H^{\tau,q}(D),
\label{8.9}
\eeq
hence in $H^{\sigma,q}$-norm, for all $\sigma<\tau$.
\end{proposition}
\demo
First note that under these hypotheses, we have
\beq
u^\nu(t)\ \text{ bounded in }\ L^\infty(D)\cap H^{\tau,q}(D).
\label{8.10}
\eeq
This bound is a direct consequence of (\ref{8.4}), (\ref{4.5}), 
and the maximum principle,
which implies $\|e^{sA}f\|_{L^\infty}\le \|f\|_{L^\infty},\ s\ge 0$.
From here, (\ref{8.8}) is a consequence of the estimate
\beq
\|u\otimes v\|_{H^{\tau,q}}\le C\|u\|_{L^\infty}\|v\|_{H^{\tau,q}}
+C\|u\|_{H^{\tau,q}}\|v\|_{L^\infty};
\label{8.11}
\eeq
cf.~\cite{T}, Chapter 13, (10.52).

To proceed, we note also that the hypothesis $u_0\in L^\infty$ plus the fact 
that $e^{tA}$ is a strongly continuous semigroup on $L^p(D)$ whenever
$p<\infty$ gives
\beq
u^\nu(t)\longrightarrow u_0\ \text{ in }\ L^p\text{-norm}, \ \ \forall\,
p<\infty.
\label{8.12}
\eeq
Hence
\beq
u^\nu(t)\otimes u^\nu(t)\longrightarrow u_0\otimes u_0\ \text{ in }\
L^p\text{-norm},\ \ \forall\, p<\infty.
\label{8.13}
\eeq
The bound (\ref{8.11}) implies $\{u^\nu(t)\otimes u^\nu(t)\}$ has $\text{weak}^*$
limit points in $H^{\tau,q}(D)$ as $\nu\searrow 0$, while (\ref{8.13}) implies
any such limit point must be $u_0\otimes u_0$.  This gives (\ref{8.9}).
The $H^{\sigma,q}$-norm convergence follows from the compactness of the
inclusion $H^{\tau,q}(D)\hookrightarrow H^{\sigma,q}(D)$.
\qed

From here we can draw conclusions about the nature of the convergence of
$p^\nu(t)$ to $p_0$, which satisfies
\beq
\nabla p_0=-\nabla_{u_0}u_0=-\text{div}(u_0\otimes u_0).
\label{8.15}
\eeq
Of course, $p^\nu(t)$ and $p_0$ are defined only up to additive constants.  
We fix these by requiring
\beq
\int\limits_D p^\nu(t,x)\, dx=0=\int\limits_D p_0(x)\, dx.
\label{8.16}
\eeq
Then we obtain the following:

\begin{proposition} \label{p8.2}
In the setting of Proposition 8.1, we have
\beq
p^\nu(t)\ \text{ bounded in }\ L^\infty(D)\cap H^{\tau,q}(D),
\label{8.17}
\eeq
for $\nu\in (0,\nu_0]$, and, as $\nu\rightarrow 0$,
\beq
p^\nu(t)\longrightarrow p_0\ \text{ in }\ H^{\sigma,q}\text{-norm},\quad
\forall\, \sigma<\tau.
\label{8.18}
\eeq
\end{proposition}

By paying closer attention to the special structure of our velocity fields,
we can improve Proposition \ref{p8.2} substantially.  Recall that
\beq
u^\nu(t,x)=s^\nu(t,|x|)x^\perp,
\label{8.19}
\eeq
with a real-valued factor $s^\nu(t,|x|)$.  Now $x^\perp=x_1\pa_{x_2}-
x_2\pa_{x_1}=r\pa_\theta$, and
\beq
\nabla_{r\pa_\theta}r\pa_\theta=-x_1\pa_{x_1}-x_2\pa_{x_2}=-x,
\label{8.20}
\eeq
so, as noted in \S{\ref{sec:1}},
\beq
\aligned
\nabla_{u^\nu}u^\nu&=-s^\nu(t,|x|)^2 x \\
&=-|u^\nu(t,x)|^2 \frac{x}{r^2},
\endaligned
\label{8.21}
\eeq
and hence
\beq
r\nabla p^\nu=|u^\nu|^2 \frac{x}{r}.
\label{8.22}
\eeq
Noting that
\beq
\frac{x}{|x|}\in L^\infty(D)\cap H^{1,p}(D),\quad
\forall\, p<2,
\label{8.23}
\eeq
we obtain via arguments used in Proposition \ref{p8.1} the following conclusion:

\begin{proposition} \label{p8.3}
In the setting of Proposition 8.1, we have
\beq
r\, \nabla p^\nu(t)\ \text{ bounded in }\ L^\infty(D)\cap
H^{\tau,q}(D),
\label{8.24}
\eeq
for $\nu\in (0,\nu_0]$, and, as $\nu\rightarrow 0$,
\beq
r \nabla p^\nu(t)\longrightarrow r\nabla p_0\ \text{ in }\
H^{\sigma,q}\text{-norm},\quad \forall\, \sigma<\tau.
\label{8.25}
\eeq
\end{proposition}

The results (\ref{8.24})--(\ref{8.25}) are weak near $x=0$, but strong away from this
point.

In case $u_0\in C^\infty(\Dbar)$, one can use methods of \S{\ref{sec:7}} to
get more precise information on $p^\nu(t,x)$, from that on $u^\nu(t,x)$.
One has from (\ref{8.1}) the smooth convergence of $\nabla p^\nu(t,x)$ to
$\nabla p_0(x)$ for $|x|\le 1-c<1$.  Results that will be presented in \S{\ref{sec:11}} 
can be applied to (\ref{8.22}) to get good control
over the boundary layer behavior of $\nabla p^\nu(t,x)$ for $1-c\le |x|\le 1$.
In particular, it will follow that 
while the pressure gradient $\nabla p^\nu(t,x)$
varies noticeably over the boundary layer, of thickness $\sim \sqrt{\nu}$,
the pressure  itself $p^\nu(t,x)$ varies only slightly over this boundary
layer.

\section{Vorticity estimates and vorticity concentration}\label{sec:9}

Here we study the vorticity
\beq
\aligned
\omega^\nu(t,x)&=\rot u^\nu(t,x) \\
&=\pa_{x_1}u^\nu_2(t,x)-\pa_{x_2}u^\nu_1(t,x)
\endaligned
\label{9.1}
\eeq
of a solution to (\ref{2.1}), i.e.,
\beq
u^\nu(t)=e^{\nu tA}u_0+\Cal{S}^\nu\alpha(t),
\label{9.2}
\eeq
under the hypothesis (\ref{1.6}) of circular symmetry, which implies as in (\ref{1.11})
that
\beq
u^\nu(t,x)=s^\nu(t,|x|)x^\perp.
\label{9.3}
\eeq
As noted in \S{\ref{sec:2}}, while (\ref{9.3}) played an important role in passing from the
Navier-Stokes equation (\ref{1.1}) to the linear equation (\ref{2.1}), it has not played
a key role in much of the subsequent linear analysis.  However, it will
play a key role in this section.  Note that (\ref{9.3}) implies
\beq
\aligned
\omega^\nu(t,x)&=\varpi^\nu(t,|x|), \\
\varpi^\nu(t,r)&=\Bigl(r\frac{d}{dr}+2\Bigr) s^\nu(t,r).
\endaligned
\label{9.4}
\eeq
In particular, $\omega^\nu(t,x)$ is circularly symmetric.

Here is our first result.

\begin{proposition} \label{p9.1}
Assume $u^\nu$ has the form (\ref{9.2})--(\ref{9.3}) with
$u_0\in L^2(D)$ and $\alpha\in C^\infty_{\flat}(\RR)$.  Then $\omega^\nu=\rot u^\nu$
belongs to $C^\infty((0,\infty)\times\Dbar)$ and satisfies the following:
\beq
\pa_t \omega^\nu=\nu\Delta \omega^\nu,\quad \text{on }\ (0,\infty)\times D,
\label{9.5}
\eeq
and
\beq
n\cdot \nabla \omega^\nu(t,x)=\frac{\alpha'(t)}{2\pi\nu},\quad
\text{on }\ (0,\infty)\times\pa D,
\label{9.6}
\eeq
where $n$ is the outward unit normal to $\pa D$.  In addition,
\beq
\int\limits_D \omega^\nu(t,x)\, dx=\alpha(t),\quad \forall\, t>0.
\label{9.7}
\eeq
\end{proposition}
\demo
Standard regularity results yield $\omega^\nu\in C^\infty
((0,\infty)\times\Dbar)$, and (\ref{9.5}) is also obvious.  Since $\omega^\nu=
-\dv (u^\nu)^\perp$, the divergence theorem gives
\beq
\int\limits_D \omega^\nu(t,x)\, dx=-\int\limits_{\pa D} n\cdot
u^\nu(t,x)^\perp\, ds=\alpha(t),
\label{9.8}
\eeq
since $u^\nu(t,x)^\perp=-\alpha(t)x/2\pi$ on $\pa D$.  Next, using (\ref{9.5}) and the
divergence theorem,
\beq
\aligned
\frac{d}{dt} \int\limits_D \omega^\nu(t,x)\, dx
&=\nu\int\limits_D \Delta\omega^\nu(t,x)\, dx \\
&=\nu\int\limits_{\pa D} x\cdot\nabla \omega^\nu(t,x)\, dx,
\endaligned
\label{9.9}
\eeq
which is $2\pi\nu$ times the left side of (\ref{9.6}), by circular symmetry.
Since the left side of (\ref{9.9}) equals $\alpha'(t)$, (\ref{9.6}) is proven.
\qed

From this follows:

\begin{proposition} \label{p9.2}
Assume
\beq
u_0\in H^1_0(D),\quad u_0(x)=s_0(|x|)x^\perp,
\label{9.10}
\eeq
and set
\beq
\omega_0=\rot u_0,\quad \omega^\nu(t)=\rot e^{\nu tA}u_0.
\label{9.11}
\eeq
Then
\beq
\omega^\nu(t)=e^{\nu tA_N}\omega_0,
\label{9.12}
\eeq
where $A_N$ is the self adjoint operator on (scalar functions in) $L^2(D)$
given by
\beq
\aligned
&\Cal{D}(A_N)=\{\omega\in H^2(D):n\cdot \nabla\omega|_{\pa D}=0\}, \\
&A_N\omega=\Delta\omega\ \text{ for }\ \omega\in \Cal{D}(A_N).
\endaligned
\label{9.13}
\eeq
\end{proposition}
\demo
Proposition \ref{p9.1} applies with $\alpha\equiv 0$, so we have
$\omega^\nu\in C^\infty((0,\infty)\times\Dbar)$ satisfying
\beq
\alignedat 2
\pa_t\omega^\nu&=\nu\Delta\omega^\nu\quad &\text{on }\ (0,\infty)\times D, \\
n\cdot\nabla\omega^\nu&=0\quad &\text{on }\ (0,\infty)\times\pa D,
\endalignedat
\label{9.14}
\eeq
hence, for $0<s<t<\infty$,
\beq
\omega^\nu(t)=e^{\nu(t-s)A_N}\omega^\nu(s).
\label{9.15}
\eeq
Also, the hypothesis $u_0\in H^1_0(D)$ implies $e^{\nu tA}u_0\in
C([0,\infty),H^1_0(D))$, since $H^1_0(D)=\Cal{D}((-A)^{1/2})$.
It follows that
\beq
\omega^\nu\in C([0,\infty),L^2(D)),\quad \omega^\nu(0)=\omega_0,
\label{9.16}
\eeq
and hence (\ref{9.12}) follows from (\ref{9.15}) in the limit $s\searrow 0$.
\qed

Generally, given
\beq
u_0\in H^1(D),\quad u_0(x)=s_0(|x|)x^\perp,
\label{9.17}
\eeq
we can write
\beq
u_0=u_{00}+af_1,\quad u_{00}\in H^1_0(D),
\label{9.18}
\eeq
with $f_1$ as in (\ref{2.19}), i.e., $f_1(x)=x^\perp/2\pi$.  In this case,
Proposition \ref{p9.2} is applicable to $e^{\nu tA}u_{00}$, but not to
\beq
\omega^\nu(t)=\rot e^{\nu tA}f_1.
\label{9.19}
\eeq
We still have (\ref{9.14})--(\ref{9.15}), but (\ref{9.16}) fails, and so does 
(\ref{9.12}).  In fact, in this case
\beq
\omega_0=\rot f_1=\frac{1}{\pi},\quad \text{so }\
e^{\nu tA_N}\omega_0\equiv \frac{1}{\pi},
\label{9.20}
\eeq
but (\ref{9.7}) holds with $\alpha\equiv 0$, so for $\omega^\nu(t)$ given by
(\ref{9.19}),
\beq
\int\limits_D \omega^\nu(t,x)\, dx=0,\quad \forall\, t>0.
\label{9.21}
\eeq
In this context, we note that Proposition \ref{p7.1} implies that, for each
compact $\Ombar\subset D$,
\beq
e^{\nu tA}f_1\longrightarrow f_1\ \text{ in }\ C^\infty(\Ombar),
\label{9.22}
\eeq
as $\nu\searrow 0$ (or as $t\searrow 0$), and hence, with $\omega^\nu,\
\omega_0$ as in (\ref{9.19})--(\ref{9.20}),
\beq
\omega^\nu(t)\longrightarrow \omega_0\ \text{ in }\ C^\infty(\Ombar)
\label{9.23}
\eeq
as $\nu\searrow 0$ (or as $t\searrow 0$).  This implies a ``concentration
phenomenon'' for the vorticity $\omega^\nu(t)$, in $\pa D$, which we discuss
in a more general context below.

Regarding the formula (\ref{9.12}), it is standard that $\{e^{sA_N}:s\ge 0\}$
is a strongly continuous contraction semigroup on $L^p(D)$, for each
$p\in [1,\infty)$.  Consequently Proposition \ref{p9.2} yields:
\beq
u_0\ \text{ as in (\ref{9.10})} \Longrightarrow \|\rot e^{\nu tA}u_0\|_{L^1(D)}
\le \|\rot u_0\|_{L^1(D)}.
\label{9.24}
\eeq
The following is a useful complement.

\begin{proposition} \label{p9.3}
Given
\beq
\alpha\in C^\infty_{\flat}(\RR),\quad u^\nu(t)=\Cal{S}^\nu\alpha(t),\quad
\omega^\nu(t)=\rot u^\nu(t),
\label{9.25}
\eeq
we have
\beq
\|\omega^\nu(t)\|_{L^1(D)}\le \|\alpha\|_{\BV([0,t])}.
\label{9.26}
\eeq
\end{proposition}
\demo
Under these hypotheses, $u^\nu,\omega^\nu\in C^\infty(\RR
\times\Dbar)$.  Take $\varphi_\ep(y)=(\ep^2+y^2)^{1/2}$, which approximates
$|y|$ as $\ep\searrow 0$.  Multiply the equation (\ref{9.5}) by
$\varphi'_\ep(\omega^\nu)$ and integrate over $D$ to obtain
\beq
\aligned
&\frac{d}{dt}\int\limits_D \varphi_\ep(\omega^\nu(t,x))\, dx \\
&=\nu\int\limits_D \Bigl[\Delta\varphi_\ep(\omega^\nu(t,x))-\varphi''_\ep
(\omega^\nu(t,x))|\nabla \omega^\nu(t,x)|^2\Bigr]\, dx \\
&\le \nu \int\limits_D \Delta \varphi_\ep(\omega^\nu(t,x))\, dx,
\endaligned
\label{9.27}
\eeq
since $\varphi''\ge 0$.  Then the divergence theorem gives
\beq
\aligned
\nu\int\limits_D \Delta\varphi_\ep(\omega^\nu(t,x))\, dx
&=\nu\int\limits_{\pa D} x\cdot\nabla\varphi_\ep(\omega^\nu(t,x))\, ds \\
&=\frac{\alpha'(t)}{2\pi} \int\limits_{\pa D} \varphi'_\ep(\omega^\nu(t,x))
\, ds,
\endaligned
\label{9.28}
\eeq
the last identity using (\ref{9.6}).  Since $|\varphi'|_\ep\le 1$, this yields
\beq
\frac{d}{dt}\int\limits_D \varphi_\ep(\omega^\nu(t,x))\, dx\le
|\alpha'(t)|.
\label{9.29}
\eeq
Consequently, for each $\ep>0$,
\beq
\int\limits_D \varphi_\ep(\omega^\nu(t,x))\, dx\le
\|\alpha\|_{\BV([0,t])}+\int\limits_D \varphi_\ep(0)\, dx.
\label{9.30}
\eeq
Taking $\ep\searrow 0$ gives (\ref{9.26}).
\qed

Returning to $e^{\nu tA}f_1$ in (\ref{9.19})--(\ref{9.20}), we note that
\beq
f_1-e^{\nu tA}f_1=\Cal{S}^\nu \chi_{\RR^+}(t),\quad \text{for }\ t>0,
\label{9.31}
\eeq
and hence, approximating $\chi_{\RR^+}$ by a sequence in $C^\infty_{\flat}(\RR)$
and passing to the limit, we get from (\ref{9.26}) that
\beq
\|\rot f_1-\rot e^{\nu tA}f_1\|_{L^1(D)}\le 1.
\label{9.32}
\eeq
By (\ref{9.20})--(\ref{9.21}),
\beq
\int\limits_D \Bigl[\rot f_1(x)-\rot e^{\nu tA}f_1(x)\Bigr]\, dx=1,
\label{9.33}
\eeq
so in fact we have identity in (\ref{9.32}), and we see the integrand in (\ref{9.33})
is $\ge 0$ on $D$, i.e.,
\beq
\rot e^{\nu tA}f_1(x)\le \frac{1}{\pi}.
\label{9.34}
\eeq
Returning to (\ref{9.32}), we have
\beq
\|\rot e^{\nu tA}f_1\|_{L^1(D)}\le 2=2\|\rot f_1\|_{L^1(D)}.
\label{9.35}
\eeq

We can put together (\ref{9.24}) and (\ref{9.35}) as follows.  Take $u_0$ as in (\ref{9.17}),
so we have (\ref{9.18}) with $a=\int_D \rot u_0(x)\, dx$, hence $\|\rot u_{00}
\|_{L^1}\le 2\|\rot u_0\|_{L^1}$.  Consequently,
\beq
u_0\ \text{ as in (\ref{9.17})} \Longrightarrow \|\rot e^{\nu tA}u_0\|_{L^1(D)}
\le 4\|\rot u_0\|_{L^1(D)}.
\label{9.36}
\eeq

We can extend the scope of (\ref{9.36}) as follows.  Set
\beq
R^1(D)=\{u\in L^2(D):u(x)=s(|x|)x^\perp,\ \rot u\in L^1(D)\}.
\label{9.37}
\eeq
An element of $R^1(D)$ is continuous on $\Dbar\setminus\{0\}$.  Set
\beq
R^1_0(D)=\bigl\{u\in R^1(D):u\bigr|_{\pa D}=0\bigr\}.
\label{9.38}
\eeq
The argument in (\ref{9.8}) readily extends to give
\beq
u\in R^1_0(D)\Longleftrightarrow u\in R^1(D)\ \text{ and }\
\int\limits_D \rot u(x)\, dx=0.
\label{9.39}
\eeq
We mention that one can apply a circularly symmetric mollifier to $\rot u$
to approximate elements of $R^1(D)$ by elements of $H^1(D)\cap R^1(D)$
and elements of $R^1_0(D)$ by elements of $H^1_0(D)\cap R^1_0(D)$.

\begin{proposition} \label{p9.4}
We have
\beq
u_0\in R^1_0(D)\Longrightarrow \|\rot e^{\nu tA}u_0\|_{L^1(D)}\le
\|\rot u_0\|_{L^1(D)},
\label{9.40}
\eeq
and
\beq
u_0\in R^1(D)\Longrightarrow \|\rot e^{\nu tA}u_0\|_{L^1(D)}\le
4\|\rot u_0\|_{L^1(D)}.
\label{9.41}
\eeq
\end{proposition}
\demo
The result (\ref{9.40}) follows from (\ref{9.24}) by a standard approximation
argument, and then (\ref{9.41}) follows by the same sort of argument as used for
(\ref{9.36}).
\qed

Here is an associated convergence result.

\begin{proposition} \label{p9.5}
We have
\beq
u_0\in R^1_0(D)\Longrightarrow \lim\limits_{\nu\searrow 0}\,
\|\rot u_0-\rot e^{\nu tA}u_0\|_{L^1(D)}=0.
\label{9.42}
\eeq
Furthermore, with $D_a=\{x\in\RR^2:|x|<a\}$, we have for each $a\in (0,1)$,
\beq
u_0\in R^1(D)\Longrightarrow \lim\limits_{\nu\searrow 0}\,
\|\rot u_0-\rot e^{\nu tA}u_0\|_{L^1(D_a)}=0.
\label{9.43}
\eeq
\end{proposition}
\demo
Take $u_0\in R^1_0(D)$.  Given $\ep>0$, there exists $v_0\in
H^1_0(D)\cap R^1_0(D)$ such that $\|\rot(u_0-v_0)\|_{L^1(D)}\le\ep$.  As in
(\ref{9.16}) we have
\beq
\rot e^{sA}v_0\in C([0,\infty),L^2(D)),
\label{9.44}
\eeq
so, using the estimate (\ref{9.40}) with $u_0$ replaced by $u_0-v_0$, we have
\beq
\limsup\limits_{\nu\searrow 0}\, \|\rot u_0-\rot e^{\nu tA}u_0\|_{L^1(D)}
\le 2\ep,
\label{9.45}
\eeq
which gives (\ref{9.42}).

More generally, given $u_0\in R^1(D)$, write $u_0=u_{00}+bf_1$, with
$u_{00}\in R^1_0(D)$.  Then (\ref{9.42}) applies to $u_{00}$, while, as we have
noted,
\beq
e^{\nu tA}f_1\longrightarrow f_1\ \text{ in }\ C^\infty(\Dbar_a),
\label{9.46}
\eeq
for each $a<1$, so (\ref{9.43}) follows.
\qed

We now delve further into the concentration phenomenon mentioned after
(\ref{9.23}).  To state it, we bring in the space of finite Borel (signed) measures
on $\Dbar$:
\beq
\Cal{M}(\Dbar)=C(\Dbar)'.
\label{9.47}
\eeq

\begin{proposition} \label{p9.6}
Given $u_0\in R^1(D)$, set
\beq
b=\int\limits_D \rot u_0(x)\, dx,
\label{9.48}
\eeq
and let $\mu$ be the rotationally invariant Borel measure on $\pa D$
of mass $1$, i.e., $1/2\pi$ times arc length on $\pa D$.  Then, for each
$t>0$,
\beq
\lim\limits_{\nu\searrow 0}\, \rot e^{\nu tA}u_0=\rot u_0-b\mu,\quad
\text{weak}^*\ \text{in } \Cal{M}(\Dbar).
\label{9.49}
\eeq
\end{proposition}
\demo
The bound (\ref{9.41}) implies $\{\rot e^{\nu tA}u_0\}$ has
$\text{weak}^*$ limit points in $\Cal{M}(\Dbar)$ as $\nu\searrow 0$.
The result (\ref{9.43}) implies any such $\text{weak}^*$ limit must be of the
form $\rot u_0-\lambda$, where $\lambda$ is a measure supported on $\pa D$.
Of course, $\lambda$ must be rotationally invariant.  Then the fact that
$\int_D \rot e^{\nu tA}u_0\, dx=0$ for each $\nu>0$ uniquely specifies
such $\lambda$ as $b\mu$.
\qed

$\text{}$ \newline
{\it Remark.}  Given $u_0\in R^1(D)$, we have
\beq
b=2\pi s_0(1).
\eeq

$\text{}$

We now pass from $\alpha\in C^\infty_{\flat}(\RR)$ to $\alpha\in \BV_{\flat}(\RR)$,
and establish the following complement to Propositions \ref{p9.5}--\ref{p9.6}.

\begin{proposition} \label{p9.7}
Assume $\alpha\in\BV_{\flat}(\RR)$ and set $v^\nu(t)=
\Cal{S}^\nu\alpha(t)$.  Then, for each $t>0$,
\beq
\|\rot v^\nu(t)\|_{L^1(D)}\le\|\alpha\|_{\BV([0,t))},
\label{9.50}
\eeq
and
\beq
\int\limits_D \rot v^\nu(t,x)\, dx=\alpha(t-).
\label{9.51}
\eeq
Furthermore, with $D_a$ as in Proposition \ref{p9.5}, we have for each $a\in(0,1)$,
\beq
\lim\limits_{\nu\searrow 0}\, \int\limits_{D_a} |\rot v^\nu(t,x)|\, dx=0.
\label{9.52}
\eeq
Therefore, with $\mu$ as in Proposition \ref{p9.6},
\beq
\lim\limits_{\nu\searrow 0}\, \rot v^\nu(t)=\alpha(t-)\mu,\quad
\text{weak}^*\ \text{in }\Cal{M}(\Dbar).
\label{9.53}
\eeq
\end{proposition}
\demo
We use the representation
\beq
v^\nu(t)=\int\limits_{[0,t)} \bigl(I-e^{\nu(t-s)A}\bigr)f_1\, d\alpha(s);
\label{9.54}
\eeq
cf.~(\ref{2.24}) and the remark following its proof.  Then (\ref{9.50}) follows from (\ref{9.32}),
(\ref{9.51}) follows from (\ref{9.33}), 
since we then have
\beq
\int\limits_D  \rot v^\nu(t,x)\, dx
=\int\limits_{[0,t)}d\alpha(s)=\alpha(t-),
\eeq
(\ref{9.52}) follows from (\ref{9.46}), and (\ref{9.53}) then
follows by the same argument as used for (\ref{9.49}).
\qed

We make the following remarks regarding the concentration of vorticity
exhibited in (\ref{9.49}) and (\ref{9.53}).  First, it implies that in considering the 
inviscid limit for the Navier-Stokes equations in domains with boundary, 
one has to deal analytically with regularity at the level of vortex sheets.
This confirms an observation made in \cite{LNP}; cf.~Remark 2 in \S{3} of that
paper.  Furthermore, for the factor $\alpha(t-)$ in (\ref{9.53}) to be nonzero,
acceleration must be applied to $\pa D$; it is this acceleration that is
responsible for the formation of the vortex sheet.

\section{Variants of results of Bona-Wu}\label{sec:10}

In \cite{BW}, J.~Bona and J.~Wu studied the small $\nu$ behavior of solutions
to (\ref{1.1})--(\ref{1.4}) in the case $\alpha\equiv 0$.  Their hypotheses on the initial
data were expressed in terms of the vorticity $\omega_0=\text{rot}\, u_0$,
which, for $u_0$ satisfying (\ref{1.4}), is radially symmetric, i.e., $\omega_0(x)
=\varpi(|x|)$.  They assumed $\varpi$ is continuous on $[0,1)$,
integrable on $[0,1]$, and satisfies $\int_0^1r\omega(r)\, dr=0$.
Under such hypotheses, it was shown that $u^\nu(t,\cdot)\rightarrow u_0$
uniformly on $\Dbar$.   Here we produce several extensions of that result.

To begin, we note that since (by (\ref{9.4}))
\beq
u_0(x)=\Bigl(\frac{1}{r^2}\int_0^r \rho\varpi(\rho)\,d\rho\Bigr)x^\perp,
\quad r=|x|,
\label{10.1}
\eeq
these hypotheses imply that
\beq
u_0\in C_*(\Dbar),
\label{10.2}
\eeq
where
\beq
C_*(\Dbar)=\bigl\{u\in C(\Dbar):u\bigr|_{\pa D}=0\bigr\}.
\label{10.3}
\eeq
Here is our first extension:

\begin{proposition} \label{p10.1}
Let $u^\nu$ satisfy (\ref{1.1})--(\ref{1.4}) with $\alpha
\equiv 0$ and (\ref{10.2}).  Given $T_0\in (0,\infty)$, we have, uniformly
in $t\in [0,T_0]$, as $\nu\rightarrow 0$,
\beq
u^\nu(t,\cdot)\longrightarrow u_0,\quad \text{uniformly on }\ \Dbar.
\label{10.4}
\eeq
\end{proposition}
\demo
Again we have (\ref{2.1}) and hence, in this situation,
\beq
u^\nu(t)=e^{\nu tA}u_0.
\label{10.5}
\eeq
The conclusion (\ref{10.4}) follows from the well known fact that $e^{sA}$ is a
strongly continuous semigroup on $C_*(\Dbar)$.
\qed

Before pursuing other results that involve the hypothesis $u_0=0$ on $\pa D$,
we present the following extension of Proposition \ref{p10.1}.

\begin{proposition} \label{p10.2}
Let $u^\nu$ satisfy (1.1)--(1.4) with $\alpha
\equiv 0$ and
\beq
u_0\in C(\Dbar).
\label{10.6}
\eeq
Given $T_0\in (0,\infty)$, we have, uniformly in $t\in [0,T_0]$, as
$\nu\rightarrow 0$,
\beq
u^\nu(t,\cdot)\longrightarrow u_0\quad \text{locally uniformly on the interior
of }\ D,
\label{10.7}
\eeq
and
\beq
\|u^\nu(t,\cdot)\|_{L^\infty(D)}\le \|u_0\|_{L^\infty(D)}.
\label{10.8}
\eeq
\end{proposition}
\demo
Again we have (\ref{10.5}).  Now $e^{sA}$ is a contraction semigroup
on $C(\Dbar)$, so we have (\ref{10.8}), but it is not strongly continuous at
$s=0$.  To get (\ref{10.7}), we argue as in the proof of Proposition \ref{p7.1}.
Let $K\subset D$ be compact.
Using a partition of unity, write
\beq
u_0=u_a+u_b,\quad u_a\in C_{b}(\Dbar),\ \ u_b=0\ \text{on a neighborhood }
U\ \text{of}\ K.
\label{10.9}
\eeq
Then $u^\nu(t)=e^{\nu tA}u_a+e^{\nu tA}u_b$, and $e^{sA}u_a\rightarrow u_a$
uniformly on $\Dbar$ as $s\searrow 0$.  It remains to show that
\beq
e^{sA}u_b\longrightarrow 0\ \text{ uniformly on }\ K, \text{ as }\ s\searrow
0.
\label{10.10}
\eeq
To see this, set
\beq
\alignedat 2
w_b(s,x)=e^{sA}&u_b(x),\quad &s\ge 0,\ x\in U, \\
&0,\quad &s<0,\ x\in U.
\endalignedat
\label{10.11}
\eeq
Then $w_b$ solves $(\pa_s-\Delta)w_b=0$ on $\RR\times U$, so hypoellipticity of
$\pa_s-\Delta$ implies
\beq
w_b\in C^\infty(\RR\times U).
\label{10.12}
\eeq
This immediately implies (\ref{10.10}) and hence (\ref{10.7}).
\qed

We can use the layer potential construction which will be carried out in 
\S{\ref{sec:11}}, especially (\ref{11.35})--(\ref{11.40}),
to sharpen Proposition \ref{p10.2} further, obtaining:

\begin{proposition} \label{p10.3}
In the setting of Proposition \ref{p10.2}, we have
\beq
u^\nu(t,x)\longrightarrow u_0(x)
\label{10.13}
\eeq
uniformly on
\beq
\{(x,\nu):|x|\le 1-\delta(\nu)\},
\label{10.14}
\eeq
whenever $\delta(\nu)$ satisfies
\beq
\frac{\delta(\nu)}{\nu^{1/2}}\longrightarrow \infty\ \text{ as }\
\nu\searrow 0.
\label{10.15}
\eeq
\end{proposition}

Let us return to the setting $u_0|_{\pa D}=0$, as hypothesized in \cite{BW}.
Another result of Theorem 2 of \cite{BW} is that, under the hypotheses made
in the first paragraph of this section, 
$\text{rot}\, u^\nu(t,\cdot)\rightarrow
\omega_0$ in $L^2(D)$-norm; equivalently, $u^\nu(t,\cdot)\rightarrow u_0$
in $H^1$-norm.  (Actually, to get such a conclusion one needs to strengthen
the hypothesis $\varpi\in L^1([0,1])$ to $\varpi\in L^2([0,1])$.)
An alternative route to such a conclusion is to note that, since
\beq
\Cal{D}((-A)^{1/2})=H^1_0(D),
\label{10.16}
\eeq
it follows that
\beq
u_0\in H^1_0(D)\Longrightarrow e^{\nu tA}u_0\rightarrow u_0\ \text{in}\
H^1\text{-norm, as }\nu\rightarrow 0.
\label{10.17}
\eeq
The following is an extension of this observation.

\begin{proposition} \label{p10.4}
Let $u_\nu$ satisfy (1.1)--(1.4) with $\alpha
\equiv 0$.  Take $\sigma\in [1,5/2)$ and assume
\beq
u_0\in H^1_0(D)\cap H^\sigma(D).
\label{10.18}
\eeq
Then, given $T_0\in (0,\infty)$, we have, uniformly in $t\in [0,T_0]$,
as $\nu\rightarrow 0$,
\beq
u^\nu(t,\cdot)\longrightarrow u_0\ \text{ in }\ H^\sigma\text{-norm.}
\label{10.19}
\eeq
\end{proposition}
\demo
This is a consequence of the fact that
\beq
\Cal{D}((-A)^{\sigma/2})=H^1_0(D)\cap H^\sigma(D),\quad
\text{for }\ \sigma\in\Bigl[1,\frac{5}{2}\Bigr),
\label{10.20}
\eeq
which, like (\ref{4.4})--(\ref{4.5}), is a special case of results of 
\cite{Se1}--\cite{Se2}.
The result (\ref{10.20}) plus what have by now become familiar arguments yields
(\ref{10.19}).
\qed

\section{Boundary layer analysis of $e^{\nu tA}u_0$}\label{sec:11}

In this section we make a detailed analysis, uniformly near $\pa D$, of
the small $\nu$ behavior of $e^{\nu tA}u_0$, uniformly on
$t\in [0,T]$, in case $u_0\in C^\infty(\Dbar)$.
This is equivalent to understanding the small $t$ behavior
of $e^{tA}u_0$.  

To start, we emphasize the case $u_0(x)=f_1(x)=x^\perp/2\pi$, later
making note of the minor modifications involved in examining the more
general case.  Note that
\beq
\alignedat 2
V(t,x)=f_1(x)-&e^{tA}f_1(x)\quad &\text{for }\ t\ge 0,\ x\in D, \\
&0\quad &\text{for }\ t<0,\ x\in D
\endalignedat
\label{11.1}
\eeq
solves
\beq
\aligned
&\pa_tV=\Delta V\quad \text{on }\ \RR\times D \\
&V\bigr|_{\RR\times\pa D}=\chi_{\RR^+}(t)f_1\bigr|_{\pa D}=g(t,x).
\endaligned
\label{11.2}
\eeq
Our task is equivalent to determining the behavior of $V(t,x)$ as $t\searrow
0$.

An argument from \S{\ref{sec:7}} is again useful here.  Namely,
the hypoellipticity of $\pa_t-\Delta$ guarantees interior regularity:
\beq
V\in C^\infty(\RR\times D).
\label{11.3}
\eeq
In particular, since $V=0$ for $t<0$, we have for each $m\in\NN,\ K\subset
D$ compact,
\beq
|V(t,x)|\le C_{m,K} t^m,\quad \text{for }\ x\in K.
\label{11.4}
\eeq
Of course, $V(t,x)=f_1(x)$ for $t>0$ and $x\in\pa D$, so there is a ``boundary
layer'' on which (\ref{11.4}) fails.

We tackle the problem of analyzing $V$ using the method of layer potentials.
Given $h$ supported in $\RR^+\times\pa D$, we set
\beq
\Cal{D}h(t,x)=\int_0^\infty \int\limits_{\pa D} h(s,y) \frac{\pa H}{\pa n_y}
(t-s,x,y)\, dS(y)\, ds,\quad t\in\RR,\ x\in D.\,
\label{11.5}
\eeq
where
\beq
H(t,x,y)=(4\pi t)^{-1} e^{-|x-y|^2/4t}\, \chi_{\RR^+}(t),
\label{11.6}
\eeq
and $n_y$ is the unit outward normal to $\pa D$ at $y$.  It is known that
\beq
\Cal{D}h\bigr|_{\RR\times\pa D}=\Bigl(\frac{1}{2}I+N\Bigr)h,
\label{11.7}
\eeq
where
\beq
Nh(t,x)=\int_0^\infty \int\limits_{\pa D} h(s,y) \frac{\pa H}{\pa n_y}
(t-s,x,y)\, dS(y)\, ds,\quad t\in\RR,\ x\in\pa D.
\label{11.8}
\eeq
Cf.~\cite{T}, Chapter 7, (13.50)--(13.55). 
Then $V$ in (11.2) is given by
\beq
V=\Cal{D}h,
\label{11.9}
\eeq
where $h$ solves
\beq
\Bigl(\frac{1}{2}I+N\Bigr)h=g,
\label{11.10}
\eeq
with $g$ given in (\ref{11.2}).  Such a solution $h$ is in fact given by
\beq
h=2(I-2N+4N^2-\cdots)g.
\label{11.11}
\eeq
Not only is this convergent (at least for small $t$), but
\beq
N\in OPS^{-1/2}_{1/2,0}(\RR\times\pa D),
\label{11.12}
\eeq
so the various terms in the series are progressively smoother.  
In fact $N$ has further structure as a singular integral operator, 
exposed in \cite{FJ} and \cite{FR}, which implies it has order
$-1/2$ on $L^p$-Sobolev spaces.

To see in an elementary manner that $N$ is a weakly singular integral 
operator, we note that for $y\in\pa D, x\in\Dbar, t>0$,
\beq
\pa_{n_y}H(t,x,y)=\frac{1}{\pi}\, \frac{1}{(4t)^2}\, n(y)\cdot (x-y)
e^{-|x-y|^2/4t}.
\label{11.13}
\eeq
This has a relatively weak singularity on $\RR^+\times\pa D\times\pa D$,
which when $D$ is the disk can be given rather explicitly, using the fact
that
\beq
n(y)=y
\label{11.14}
\eeq
for $y\in\pa D$, so $n(y)\cdot(x-y)=x\cdot y-1$.  Also $|x-y|^2=2-2x\cdot y$,
so $n(y)\cdot(x-y)=-|x-y|^2/2$, for $x,y\in\pa D$.  Hence
\beq
\aligned
\pa_{n_y}H(t,x,y)&=\frac{1}{8\pi t}\, \frac{|x-y|^2}{4t}\,
e^{-|x-y|^2/4t} \\ &=\frac{1}{8\pi t}\, \Phi
\Bigl(\frac{|x-y|}{\sqrt{4t}}\Bigr),\quad
\text{for }\ x,y\in\pa D,\ t>0,
\endaligned
\label{11.15}
\eeq
where
\beq
\Phi(\lambda)=\lambda^2 e^{-\lambda^2}.
\label{11.16}
\eeq
Note that this is less singular than its counterpart where $y\in\pa D$
and $x\in D$ approaches $y$ radially, by a factor of $|x-y|$.  From here
on, we denote by $N(t,x,y)$ the function on $\RR\times\Dbar\times\pa D$
given by (\ref{11.13}) for $t>0$, and vanishing for $t<0$.

Clearly we have
\beq
\|N(t,x,\cdot)\|_{L^1(\pa D)}\le Ct^{-1/2},\quad t\in\RR^+,\ x\in\pa D.
\label{11.17}
\eeq
Hence, with $g$ as in (11.2), $t>0$,
\beq
|Ng(t,x)|\le Ct^{1/2}\|g\|_{L^\infty}.
\label{11.18}
\eeq
By contrast,
\beq
\|N(t,x,\cdot)\|_{L^1(\pa D)}\le Ct^{-1},\quad t\in\RR^+,\ x\in \Dbar.
\label{11.19}
\eeq

We can deduce that the solution $h$ to (\ref{11.10}) satisfies
\beq
h=2g+h^b,
\label{11.20}
\eeq
with $h^b(t,x)$ supported in $t\in\RR^+$ and (at least for small $t$)
\beq
|h^b(t,x)|\le Ct^{1/2}.
\label{11.21}
\eeq
Hence
\beq
V(t,x)=\Cal{D}h(t,x)=2\Cal{D}g(t,x)+\Cal{D}h^b(t,x).
\label{11.22}
\eeq
To estimate $\Cal{D}h^b(t,x)$, (\ref{11.19}) is not so useful; instead we
argue as follows.  Denote by $\PI$ the solution operator to (\ref{11.2}), so
(\ref{11.7}) gives
\beq
\PI g=\Cal{D}h,\quad g=\Bigl(\frac{1}{2}I+N\Bigr)h.
\label{11.23}
\eeq
Similarly
\beq
\Cal{D}h^b=\PI g^b,\quad g^b=\Bigl(\frac{1}{2}I+N\Bigr)h^b.
\label{11.24}
\eeq
As usual, $g^b(t,x)$ is supported in $t\in\RR^+$.
Also, (\ref{11.18}) and (\ref{11.21}) give
\beq
|g^b(t,x)|\le Ct^{1/2},
\label{11.25}
\eeq
for $t>0$.  Then the maximum principle for solutions to the heat equation
gives
\beq
|\PI g^b(t,x)|\le C t^{1/2},
\label{11.26}
\eeq
and by (\ref{11.24}) this is the estimate we have on $\Cal{D}h^b(t,x)$.
We have established the following.

\begin{proposition} \label{p11.1}
For $V$ given by (\ref{11.1}) and $g$ in (\ref{11.2}), we have
\beq
|V(t,x)-2\Cal{D}g(t,x)|\le Ct^{1/2},\quad \forall\, x\in\Dbar.
\label{11.27}
\eeq
\end{proposition}

To get finer approximations to $V(t,x)$, we use more terms in (\ref{11.11}).  Write
\beq
h=2g_k+h^b_k,
\label{11.28}
\eeq
where
\beq
g_k=\sum\limits_{j=0}^k (-2N)^jg,\quad h^b_k=2(-2N)^{k+1}
\sum\limits_{j=0}^\infty (-2N)^j g.
\label{11.29}
\eeq
We have, for small $t>0$,
\beq
\aligned
|N^jg(t,x)|&\le (Ct^{1/2})^j \|g\|_{L^\infty}, \\
|h^b_k(t,x)|&\le (Ct^{1/2})^{k+1} \|g\|_{L^\infty}.
\endaligned
\label{11.30}
\eeq
Furthermore, $N^j$ has smoothing properties, leading to the fact that
$N^jg$ and $h^b_k$ are supported in $t\in\RR^+$ and
\beq
N^jg\in C^{j/2-\ep}(\RR\times\pa D),\quad h^b_k\in C^{(k+1)/2-\ep}
(\RR\times\pa D).
\label{11.31}
\eeq
Now, parallel to (\ref{11.22}), we have
\beq
V(t,x)=2\Cal{D}g_k(t,x)+\Cal{D}h^b_k(t,x),
\label{11.32}
\eeq
where $\Cal{D}h^b_k(t,x)$ is supported on $\RR^+\times\Dbar$ and, parallel to
(\ref{11.24}),
\beq
\Cal{D}h^b_k=\PI g^b_k,\quad g^b_k=\Bigl(\frac{1}{2}I+N\Bigr)h^b_k.
\label{11.33}
\eeq
Thus $g^b_k\in C^{(k+1)/2-\ep}(\RR\times\pa D)$, and consequently $\PI g^b_k$
has this degree of regularity on $\RR\times\Dbar$.  In conclusion:

\begin{proposition} \label{p11.2}
In the setting of Proposition \ref{p11.1}, define $g_k$ by (\ref{11.29}).
Then, for each $k\in\NN$, we have (\ref{11.32}), with 
remainder
\beq
\Cal{D}h^b_k\in C^{(k+1)/2-\ep}(\RR\times\Dbar),
\label{11.34}
\eeq
supported in $\RR^+$.
\end{proposition}

Of course there is a similar treatment of $e^{tA}u_0$ when $u_0\in C^\infty
(\Dbar)$.  One can take an extension $u_0\in C^\infty(\RR^2)$
(polynomially bounded, say), set
\beq
U_0(t,x)=e^{t\Delta}u_0(x),\quad \text{on }\ \RR^+\times\RR^2,
\label{11.35}
\eeq
and then note that
\beq
\alignedat 2
W(t,x)=U_0(t,x)&-e^{tA}u_0(x)\quad &\text{for }\ t\ge 0,\ x\in D \\
&0\qquad &\text{for }\ t<0,\ x\in D
\endalignedat
\label{11.36}
\eeq
solves
\beq
\aligned
&\pa_tW=\Delta W\quad \text{on }\ \RR\times D, \\
&W\bigr|_{\RR\times\pa D}=\chi_{\RR^+}(t) U_0(t,x)\bigr|_{\pa D}=
\tilde{g}(t,x),
\endaligned
\label{11.37}
\eeq
from which one has straightforward parallels of (\ref{11.3})--(\ref{11.34}).  One has
(for $t$ small)
\beq
W=\Cal{D}\tilde{h},
\label{11.38}
\eeq
where $\tilde{h}$ solves
\beq
\Bigl(\frac{1}{2}I+N\Bigr)\tilde{h}=\tilde{g},
\label{11.39}
\eeq
so
\beq
\tilde{h}=2(I-2N+4N^2-\cdots)\tilde{g}.
\label{11.40}
\eeq
As with $g$, we have that $\tilde{g}$ is piecewise smooth, with a jump across
$\{t=0\}$.

We can go further, and make the construction (\ref{11.35})--(\ref{11.40}) for more general
$u_0$, such as $u_0\in C(\Dbar)$.  In this case, one takes a polynomially bounded 
extension $u_0\in C(\RR^2)$ to define $U_0$ in (\ref{11.35}).  Then $\tilde{g}$ in 
(\ref{11.37}) is piecewise continuous, with a jump across $\{t=0\}$.  With $\tilde{h}$
given by (\ref{11.38})--(\ref{11.40}), estimates on $\tilde{h}$ parallel to those
on $h$ (given by (\ref{11.11})) hold.  In particular, parallel to (\ref{11.20}) we have 
$\tilde{h}=2\tilde{g}+\tilde{h}^b$, and $\tilde{h}^b$ has a treatment analogous to 
(\ref{11.21})--(\ref{11.26}).  Thus we have the following analogue of Proposition \ref{p11.1},
which we state explicitly, since it implies Proposition \ref{p10.3}, as advertized in 
\S{\ref{sec:10}}.

\begin{proposition} \label{p11.3}
Given $u_0\in C(\Dbar)$, define $W(t,x)$ by (11.35)--(11.36), and define $\tilde{g}(t,x)$
by (11.37).  Then we have
\beq
|W(t,x)-2\Cal{D}\tilde{g}(t,x)|\le C t^{1/2},\quad \forall\, x\in\Dbar.
\label{11.41}
\eeq
\end{proposition}

Similarly Proposition \ref{p11.2} extends to this setting, with minimal changes 
in the proof.

\section{Boundary layer analysis of $\Cal{S}^\nu\alpha(t)$}\label{sec:12}

Let us assume $\alpha\in\BV_{\flat}(\RR)$.
We can apply (\ref{11.22})--(\ref{11.27}) to analyze
\beq
\aligned
\Cal{S}^\nu\alpha(t)&=\int_0^t (I-e^{\nu(t-s)A})f_1\, d\alpha(s) \\
&=\int_0^t V(\nu(t-s),x)\, d\alpha(s),
\endaligned
\label{12.1}
\eeq
obtaining
\beq
\Cal{S}^\nu\alpha(t)=2\int_0^t \Cal{D}g(\nu(t-s),x)\, d\alpha(s)
+\Cal{R}_\nu(t),
\label{12.2}
\eeq
with $g$ given by (\ref{11.2}) and
\beq
\|\Cal{R}_\nu(t)\|_{L^\infty(D)}\le Ct^{1/2}\nu^{1/2}\|\alpha\|_{\TV([0,t])}.
\label{12.3}
\eeq

Next, we can apply (\ref{11.32})--(\ref{11.34}) to analyze
\beq
\aligned
\Cal{S}^\nu\alpha(t)&=-\int_0^t \nu Ae^{\nu(t-s)A}f_1\, \alpha(s)\, ds \\
&=-\nu\Delta f_1\int_0^t \alpha(s)\, ds+\nu\int_0^t \Delta V(\nu(t-s),x)
\alpha(s)\, ds.
\endaligned
\label{12.4}
\eeq
Note that $\Delta f_1=0$, so we have
\beq
\Cal{S}^\nu\alpha(t)=2\nu \int_0^t \Delta 
\Cal{D}g_k(\nu(t-s),x)\alpha(s)\, ds
+\nu \Cal{R}^2_{\nu,k}(t,x),
\label{12.5}
\eeq
with
\beq
\Cal{R}^2_{\nu,k}\in C^{(k+1)/2-2-\ep}(\RR\times\Dbar),
\label{12.6}
\eeq
supported in $t\in\RR^+$, given $\alpha\in L^1_{b}(\RR)$.  In
particular,
\beq
\|\nu\Cal{R}^2_{\nu,k}(t,\cdot)\|_{L^\infty(D)}\le c\nu^{(k+1)/2-1-\ep}
\|\alpha\|_{L^1([0,t])}.
\label{12.7}
\eeq

The significance of these estimates is that the principal term on the right 
side of (\ref{12.2}) is an explicit integral (and its counterpart in
(\ref{12.5}) is more or less explicit).

\section{Concentric rotating circles} \label{sec:13}

Here we extend the analysis of the previous sections to the case where
the disk $D$ is replaced by the annulus
\beq
\Cal{A}=\{x\in\RR^2:\rho<|x|<1\},
\label{13.1}
\eeq
for some $\rho\in (0,1)$.  
The results in this section are an extension of work described in \cite{FGP}.

Thus we consider solutions to the Navier-Stokes
equations (\ref{1.1}) on $\RR^+\times\Cal{A}$, with no-slip boundary data on the
two components of $\pa\Cal{A}$, which might be rotating independently:
\beq
\aligned
u^\nu(t,x)=\ &\frac{\alpha_1(t)}{2\pi} x^\perp,\quad |x|=1,\ t>0, \\
&\frac{\alpha_2(t)}{2\pi} x^\perp,\quad |x|=\rho,\ t>0,
\endaligned
\label{13.2}
\eeq
and with circularly symmetric initial data
\beq
u^\nu(0,x)=u_0(x),\quad \dv u_0=0,\quad u_0\, \| \, \pa\Cal{A},
\label{13.3}
\eeq
where again circular symmetry is defined by (\ref{1.4}), and entails the formulas
(\ref{1.6}) and (\ref{1.7}).  Proposition \ref{p1.1} immediately extends to this setting, so
$u^\nu(t,x)$ satisfies (\ref{1.11}) for all $t>0$ and is specified as the solution
to the linear PDE
\beq
\pa_t u^\nu=\nu\Delta u^\nu,
\label{13.4}
\eeq
on $\RR^+\times\Cal{A}$, with boundary data (\ref{13.2}) and initial data
(\ref{13.3}).  The material of \S{2} easily extends.
We can represent the solution to (\ref{13.2})--(\ref{13.4}) as
\beq
u^\nu(t)=e^{\nu t\Delta}u_0+\Cal{S}^\nu(\alpha_1,\alpha_2)(t),
\label{13.5}
\eeq
where $A$ is given by
\beq
\Cal{D}(A)=H^2(\Cal{A})\cap H^1_0(\Cal{A}),\quad Au=\Delta u\ \text{ for }\
u\in\Cal{D}(A),
\label{13.6}
\eeq
and
\beq
\Cal{S}^\nu:C^\infty_\flat(\RR)\oplus C^\infty_\flat(\RR)\longrightarrow
C^\infty_\flat(\RR\times\overline{\Cal{A}})
\label{13.7}
\eeq
is defined as $\Cal{S}^\nu(\alpha_1,\alpha_2)=v^\nu$, where $v^\nu$ is the
solution to (\ref{13.2})--(\ref{13.4}) that vanishes for $t<0$ (as in
(\ref{2.5})).  As before, we have extensions of $\Cal{S}^\nu$ such as
\beq
\Cal{S}^\nu:C_\flat(\RR)\oplus C_\flat(\RR)\longrightarrow
C_\flat(\RR\times\overline{\Cal{A}}),
\label{13.8}
\eeq
and also analogues of (\ref{2.10})--(\ref{2.15}).  We obtain analogues of
(\ref{2.17})--(\ref{2.22}) as follows.
With $v^\nu=\Cal{S}^\nu(\alpha_1,\alpha_2)$ defined above, set
\beq
w^\nu(t,x)=v^\nu(t,x)-\Phi(t,x)
\label{13.9}
\eeq
on $[0,\infty)\times\Cal{A}$, where $\Phi(t,x)$ is defined for each $t\ge 0$
by
\beq
\Delta\Phi(t,\cdot)=0\ \text{ on }\ \Cal{A},\quad \Phi(t,x)=
\frac{\alpha_j(t)}{2\pi} x^\perp\ \text{ on }\ B_j,
\label{13.10}
\eeq
where $B_1=\{x\in\RR^2:|x|=1\},\ B_2=\{x\in\RR^2:|x|=\rho\}$.  Then $w^\nu$
solves
\beq
\pa_t w^\nu=\nu\Delta w^\nu-\pa_t\Phi,\quad w^\nu(0,x)=0,\quad
w^\nu\bigr|_{\RR^+\times\pa\Cal{A}}=0,
\label{13.11}
\eeq
so by Duhamel's formula we have the following variant of (\ref{2.20}):
\beq
w^\nu(t)=-\int_0^t e^{\nu(t-s)A} \pa_s\Phi(s,x)\, ds.
\label{13.12}
\eeq
Hence
\beq
\aligned
\Cal{S}^\nu(\alpha_1,\alpha_2)(t)&=\Phi(s,x)-\int_0^t e^{\nu(t-s)A}
\pa_s\Phi(s,x)\, ds \\
&=\int_0^t (I-e^{\nu(t-s)A})\pa_s \Phi(s,x)\, ds.
\endaligned
\label{13.13}
\eeq
As in \S{\ref{sec:2}}, we first get these identities for $\alpha_j\in C^\infty_\flat
(\RR)$, and then we can extend the validity of these formulas via limiting
arguments.  We can obtain formulas more closely resembling (\ref{2.21}) as follows.
Vector fields on $\Cal{A}$ of the form $s_0(|x|)x^\perp$ that are harmonic
are linear combinations of $x^\perp$ and $|x|^{-2}x^\perp$, so $\Phi(t,x)$,
defined by (\ref{13.5}), is given by
\beq
\Phi(t,x)=\beta_1(t)f_1(x)+\beta_2(t)f_2(x),
\label{13.14}
\eeq
with
\beq
f_1(x)=\frac{x^\perp}{2\pi},\quad f_2(x)=\frac{x^\perp}{2\pi |x|^2},
\label{13.15}
\eeq
and $\beta_j(t)$ given by
\beq
\aligned
\beta_1(t)+\beta_2(t)&=\alpha_1(t), \\
\beta_1(t)+\frac{\beta_2(t)}{\rho^2}&=\alpha_2(t).
\endaligned
\label{13.16}
\eeq
Solving for $\beta_j$ and plugging into (\ref{13.13}), we obtain
\beq
\aligned
&\Cal{S}^\nu(\alpha_1,\alpha_2)(t) \\
&=\int_0^t \Bigl(I-e^{\nu(t-s)A}\Bigr)\Bigl[\frac{f_1-\rho^2 f_2}{1-\rho^2}
\alpha'_1(s)-\frac{\rho^2(f_1-f_2)}{1-\rho^2}\alpha'_2(s)\Bigr]\, ds,
\endaligned
\label{13.17}
\eeq
first for $\alpha_j\in C^\infty_\flat(\RR)$, then in more general cases.
For example, parallel to Proposition \ref{p2.1}, we have
\beq
\Cal{S}^\nu:{\BV}_\flat(\RR)\oplus \BV_\flat(\RR)\longrightarrow
C_\flat(\RR,X),
\label{13.18}
\eeq
whenever $X$ is a Banach space of functions on $\Cal{A}$ such that $f_1, f_2
\in X$ and $\{e^{tA}:t\ge 0\}$ is a strongly continuous semigroup on $X$.
In such a case,
\beq
\aligned
&\Cal{S}^\nu(\alpha_1,\alpha_2)(t) \\
&=\int\limits_{I(t)} \Bigl(I-e^{\nu(t-s)A}\Bigr)\Bigl[
\frac{f_1-\rho^2 f_2}{1-\rho^2}\, d\alpha_1(s)
-\frac{\rho^2(f_1-f_2)}{1-\rho_2}\, d\alpha_2(s)\Bigr],
\endaligned
\label{13.19}
\eeq
where we can take $I(t)=[0,t]$ or $I(t)=[0,t)$.  Results of
\S\S{\ref{sec:3}--\ref{sec:8}}
extend in a straightforward way to the current setting.

To extend the results of \S{9}, we need to do a little more work.
First we have the following variant of Proposition \ref{p9.1}.

\begin{proposition} \label{p13.1}
Assume $u^\nu(t,x)=s^\nu(t,|x|)x^\perp$ has
the form (\ref{13.5}) with $u_0\in L^2(\Cal{A})$ and $\alpha_j\in C^\infty_\flat
(\RR)$.  Then $\omega^\nu=\rot u^\nu$ belongs to $C^\infty((0,\infty)
\times\overline{\Cal{A}})$ and satisfies the following:
\beq
\pa_t\omega^\nu=\nu\Delta \omega^\nu\ \text{ on }\ (0,\infty)\times\Cal{A},
\label{13.20}
\eeq
and
\beq
n\cdot\nabla\omega^\nu(t,x)=(-1)^{j-1} |x| \frac{\alpha'_j(t)}{2\pi\nu}
\ \text{ on }\ B_j,
\label{13.21}
\eeq
with $B_j$ as in (13.10) and $n$ the outward unit normal to $\pa\Cal{A}$.
Also
\beq
\int\limits_{\Cal{A}} \omega^\nu(t,x)\, dx=\alpha_1(t)-\rho\alpha_2(t).
\label{13.22}
\eeq
\end{proposition}

\demo
The results (\ref{13.20}) and (\ref{13.22}) are proven just as in
Proposition \ref{p9.1}.  However, (\ref{13.21}) does not follow as easily as
(\ref{9.6}), because $\pa\Cal{A}$ has two components.  Instead, we calculate
as follows.  With $\tilde{n}=x/|x|$, we have
\beq
\aligned
\tilde{n}\cdot\nabla\omega^\nu(t,x)&=\frac{\pa}{\pa r}\varpi^\nu(t,r) \\
&=\frac{\pa}{\pa r} \Bigl(r\frac{\pa}{\pa r}+2\Bigr)s^\nu(t,r) \\
&=\Bigl(r\frac{\pa^2}{\pa r^2}+3\frac{\pa}{\pa r}\Bigr)s^\nu(t,r).
\endaligned
\label{13.23}
\eeq
Note that $\tilde{n}=(-1)^{j-1}n$ on $B_j$.  Also we have
\beq
\aligned
\Delta u^\nu&=(\Delta s^\nu)x^\perp+2\nabla s^\nu\cdot\nabla x^\perp \\
&=\Bigl(\Delta s^\nu +\frac{2}{r}\pa_r s^\nu\Bigr)x^\perp \\
&=\Bigl(\frac{\pa^2}{\pa r^2}+\frac{3}{r}\frac{\pa}{\pa r}\Bigr)s^\nu\
x^\perp,
\endaligned
\label{13.24}
\eeq
while
\beq
\aligned
\Delta u^\nu\bigr|_{\RR^+\times B_j}&=\frac{1}{\nu} \pa_tu^\nu
\bigr|_{\RR^+\times B_j} \\
&=\frac{\alpha'_j(t)}{2\pi\nu}\, x^\perp.
\endaligned
\label{13.25}
\eeq
Comparison of (\ref{13.23})--(\ref{13.25}) yields (\ref{13.21}).
\qed

From here we readily extend Proposition \ref{p9.2}, obtaining
\beq
\aligned
u_0\in H^1_0(\Cal{A}),\quad &u_0(x)=s_0(|x|)x^\perp,\quad \omega_0=\rot u_0,
\\ &\omega^\nu(t)=\rot e^{\nu tA}u_0\Longrightarrow \omega^\nu(t)=
e^{\nu tA_N}\omega_0.
\endaligned
\label{13.26}
\eeq
and the results described in Propositions \ref{p9.3}--\ref{p9.5} readily
extend to the current setting.  In particular, with
\beq
R^1(\Cal{A})=\{u_0\in L^2(\Cal{A}):u_0(x)=s_0(|x|)x^\perp,\, \rot u_0
\in L^1(\Cal{A})\},
\label{13.27}
\eeq
we have
\beq
u_0\in R^1(\Cal{A})\Longrightarrow \|\rot e^{\nu tA}u_0\|_{L^1(\Cal{A})}
\le 4 \|\rot u_0\|_{L^1(\Cal{A})}.
\label{13.28}
\eeq
Consequently, given $u_0\in R^1(\Cal{A})$, the family $\{\rot e^{\nu tA}u_0\}$
has $\text{weak}^*$ limit points in $\Cal{M}(\overline{\Cal{A}})$ as
$\nu\searrow 0$, for each $t\in (0,\infty)$.  If we compare the integral
\beq
\int\limits_{\Cal{A}} \rot u_0(x)\, dx=2\pi\bigl[s_0(1)-\rho s_0(\rho)\bigr]
\label{13.29}
\eeq
which need not be zero, with
\beq
\int\limits_{\Cal{A}} \rot e^{\nu tA}u_0\, dx=0,\quad \forall\, \nu,t>0,
\label{13.30}
\eeq
which follows from (\ref{13.22}), we see there is a concentration phenomenon,
such as described in Propositions \ref{p9.6}--\ref{p9.7}.  We will establish
the following variant of Proposition \ref{p9.6}.

\begin{proposition} \label{p13.2}
Let $\mu_j$ be the rotationally invariant
Borel measures of mass 1 on the components $B_j$ of $\pa\Cal{A}$.  Then,
given $u_0\in R^1(\Cal{A})$, we have for each $t>0$,
\beq
\lim\limits_{\nu\searrow 0}\, \rot e^{\nu tA}u_0=\rot u_0-2\pi s_0(1)\mu_1
+2\pi\rho s_0(\rho)\mu_2,
\label{13.31}
\eeq
$\text{weak}^*$ in $\Cal{M}(\overline{\Cal{A}})$.
\end{proposition}

We see from (\ref{13.29})--(\ref{13.30}) that the measure concentrated on $\pa\Cal{A}$
on the right side of (\ref{13.31}) has the correct integral against $1$.  The
fact that $\pa\Cal{A}$ has two connected components requires us to devote
greater effort than was needed in Proposition \ref{p9.6} to proving (\ref{13.31}).
We will prove (\ref{13.31}) with the aid of the following localization result,
which is of independent interest.  To state it, let $A_D$ stand for the 
operator denoted $A$ in \S{\ref{sec:9}}:
\beq
\Cal{D}(A_D)=H^2(D)\cap H^1_0(D),\quad A_Dv=\Delta v\ \text{ for }\
v\in\Cal{D}(A_D).
\label{13.32}
\eeq

\begin{proposition} \label{p13.3}
Consider $u_0\in R^1(\Cal{A})$ and $v_0\in R^1(D)$,
and assume
\beq
u_0(x)=v_0(x)\ \text{ for }\ x\in\Cal{O}=\{x\in\Cal{A}:|x|>(1+\rho)/2\}.
\label{13.33}
\eeq
Also set $\Cal{O}_1=\{x\in\Cal{A}:|x|>(2+\rho)/3\}$.  Then, for each $t>0$,
\beq
e^{\nu tA}u_0-e^{\nu tA_D}v_0\longrightarrow 0\ \text{ in }\
C^\infty(\overline{\Cal{O}}_1),
\label{13.34}
\eeq
as $\nu\searrow 0$.
\end{proposition}

\demo
Define $W$ on $\RR\times\Cal{O}$ by
\beq
\alignedat 2
W(t,x)=e^{tA}u_0&-e^{tA_D}v_0,\quad &t\ge 0, \\
&0,\qquad &t<0.
\endalignedat
\label{13.35}
\eeq
Then $e^{\nu tA}u_0(x)-e^{\nu tA_D}v_0(x)=W(\nu t,x)$.  Note that $W$ in
(\ref{13.35}) solves
\beq
\pa_t W=\Delta W\ \text{ on }\ \RR\times\Cal{O},\quad
W\bigr|_{\RR\times B_1}=0.
\label{13.36}
\eeq
Standard results on regularity up to the boundary give
\beq
W\in C^\infty(\RR\times\overline{\Cal{O}}_1),
\label{13.37}
\eeq
which in turn gives (\ref{13.34}).
\qed

We use Proposition \ref{p13.3} to prove Proposition \ref{p13.2}.  Any $\text{weak}^*$
limit point of $\{\rot e^{\nu tA}u_0\}$ as $\nu\searrow 0$ must have the
for $\rot u_0+\lambda$, where $\lambda$ is a signed measure supported on
$\pa\Cal{A}$.  If we take $v_0$ as in Proposition \ref{p13.3} and apply Proposition
\ref{p9.6}, we have $\rot e^{\nu tA_D}v_0$ tending to $\rot v_0-2\pi s_0(1)\mu_1$
$\text{weak}^*$ in $\Cal{M}(\Dbar)$.  By (\ref{13.34}) we see that $\rot u_0+
\lambda$ must coincide with this measure when restricted to
$\overline{\Cal{O}}_1$.  Given this, (\ref{13.31}) now follows from the previous
comments about the integral against $1$, plus rotational invariance.

In a similar fashion we have the following variant of Proposition \ref{p9.7}.

\begin{proposition} \label{p13.4}
Assume $\alpha_j\in\BV_\flat(\RR)$ and set
$v^\nu(t)=\Cal{S}^\nu(\alpha_1,\alpha_2)(t)$.  Then we have, $\text{weak}^*$
in $\Cal{M}(\overline{\Cal{A}})$,
\beq
\lim\limits_{\nu\searrow 0}\, \rot v^\nu(t)=\alpha_1(t-)\mu_1-\rho
\alpha_2(t-)\mu_2,
\label{13.38}
\eeq
for each $t>0$, with $\mu_j$ as in Proposition \ref{p13.2}.
\end{proposition}

This concludes our discussion of extensions of results of \S{\ref{sec:9}}.
Extensions of results of \S\S{\ref{sec:10}--\ref{sec:12}} to the current
setting are straightforward.

\end{document}